\documentstyle[amssymb,amsfonts,12pt]{amsart}

\newcommand{\Chow}{\mbox{\it Chow}\,}
\newcommand{\DOT}{\setlength{\unitlength}{1pt}\begin{picture}(2.5,2)(1,1)
\put(1,2){\circle*{2}}\end{picture}}
\newcommand{\Fdot}{{F\!_{\DOT}}}
\newcommand{\Fpdot}{{{F\!_{\DOT}}'}}
\newcommand{\Hom}{\mbox{Hom}}
\newcommand{\Mdot}{{M_{\,\DOT}}}
\newcommand{\Span}[1]{\langle #1 \rangle}
\newcommand{\QED}{
\setlength{\unitlength}{1.0pt}%
\begin{picture}(7.5,7.5)
\put(0,-2.5){\rule{2.5pt}{2.5pt}}
\put(0,0){\rule{5pt}{2.5pt}}
\put(0,2.5){\rule{7.5pt}{2.5pt}}
\end{picture}\vspace{10pt}}

\begin{document}

\title[Pieri's formula via explicit rational equivalence]{Pieri's
formula via explicit rational equivalence}

\author{Frank Sottile}

\address{ Department of Mathematics\\ University of Toronto\\
	100 St. George Street\\ Toronto, Ontario M5S 3G3
	\\ Canada\\ (416) 978-4031} 
	\email{sottile@@math.toronto.edu} 
	\date{8 May 1996}
        \thanks{Research supported in part by NSERC grant \#
        OGP0170279} 
	\thanks{revised version of alg-geom preprint \# 9601006}
\maketitle

\begin{abstract}
Pieri's formula describes the intersection product of a Schubert 
cycle by a special Schubert cycle 
on a Grassmannian.
Classically, it is proven by calculating an appropriate 
triple intersection of Schubert varieties.
We present a new geometric proof,
exhibiting an explicit chain of rational equivalences 
from  a suitable sum of distinct Schubert varieties
to the intersection of a Schubert variety with a special
Schubert variety. 
The geometry of these rational equivalences indicate a link to a
combinatorial proof of Pieri's formula using Schensted
insertion. 
\end{abstract}

\section{Introduction}
 
Pieri's formula asserts that the product of a Schubert class and a special
Schubert class is a sum of certain other Schubert classes, each with
coefficient 1.
This determines the multiplicative structure of the Chow ring of a
Grassmann variety.
Pieri's formula also arises in algebraic combinatorics and
representation theory~\cite{Fulton_tableaux}, and has independent proofs in
each context.
One such proof~\cite{Fulton_tableaux} in combinatorics involves Schensted
insertion~\cite{Schensted}.
Its geometric proof (for example, in Hodge and
Pedoe~\cite{Hodge_Pedoe}) 
involves studying an intersection of three Schubert varieties and
invoking Poincar{\'e} duality to obtain the desired sum. 
We present a new geometric proof of this formula,
explicitly describing a sequence of deformations (inducing
rational equivalence) that transform a general intersection of a
Schubert variety with a special Schubert variety into a sum of
distinct Schubert varieties.
The geometry of these deformations is quite interesting and their
form parallels the combinatorial proof of Pieri's formula using
Schensted insertion.

Let ${\bf G}_mV$ be the Grassmannian of  
$m$-dimensional subspaces of an $n$-dimen\-sional vector
space $V$. 
A decreasing sequence $\alpha$ of length $m$, 
($n\geq \alpha_1>\cdots>\alpha_m\geq 1$),
and a complete flag $\Fdot$ in $V$  together determine  a
Schubert subvariety  $\Omega_\alpha\Fdot$ of ${\bf G}_mV$.
Special Schubert varieties $\Omega_L$ are those Schubert
varieties given by the single condition that an $m$-plane
intersect a given linear subspace $L$ non-trivially.
For any subscheme $X$ of  ${\bf G}_mV$, let $[X]$ be the cycle
class of $X$ in the Chow ring of ${\bf G}_mV$.
Pieri's formula asserts
\begin{equation}\label{eq:pieri}
[\Omega_\alpha\Fdot]\cdot[\Omega_L]\  =\  
\sum [\Omega_\gamma\Fdot],
\end{equation}
the sum over all sequences $\gamma$  with 
$\gamma_1\geq \alpha_1>\gamma_2\geq 
\cdots>\gamma_m\geq \alpha_m$ where  
$\sum \gamma_i-\alpha_i$ is equal to the codimension $b$ of $\Omega_L$.
Let $\alpha *b$ denote this set of sequences
and let $Y_{\alpha,b}$ be the cycle
$\sum_{\gamma \in \alpha*b}\Omega_\gamma\Fdot$.
Let $\Chow {\bf G}_mV$ be the  Chow variety of ${\bf G}_mV$
and let ${\cal G}\subset \Chow {\bf G}_mV$ 
be the set of 
cycles $\Omega_\alpha\Fdot \bigcap \Omega_L$ for 
all $L$ of a fixed dimension such that the intersection is
generically  transverse.

Our proof involves a partial compactification of ${\cal G}$ in 
$\Chow{\bf G}_mV$ with $b+1$ rational strata, each an orbit of
the Borel subgroup of $GL(V)$ stabilizing $\Fdot$, hence consisting
of isomorphic cycles.
The 0th stratum is dense in ${\cal G}$ and cycles in the $i$th stratum have 
components 
$X_\beta$ indexed by $\beta \in \alpha* i$, where the component
$X_\beta$ is a subvariety of $\Omega_\beta\Fdot$.
Passing from one stratum to the next, each component $X_\beta$
deforms into some components of cycles in the next stratum.
The `history' of each component $\Omega_\gamma\Fdot$ of 
$Y_{\alpha,b}$ through this process gives a chain in the Bruhat
order of Schubert varieties, recording which component at each
stage gave rise to $\Omega_\gamma\Fdot$. 
This leads to the following interpretation of Pieri's formula:
The sum in (\ref{eq:pieri}) is over a certain set of chains in
the Bruhat order which begin at $\alpha$, the chain with endpoint
$\gamma$ recording the history of the cycle $\Omega_\gamma\Fdot$
in the sequence of deformations.
In \S\ref{sec:schensted}, we show how this is similar to  
a combinatorial proof of Pieri's formula based on Schensted
insertion.  

In~\cite{sottile_real_lines}, these deformations were constructed in the
special case of ${\bf G}_2V$ and applied to obtain a completely geometric 
understanding of  intersections of Schubert subvarieties of 
${\bf G}_2V$ in terms of explicit, multiplicity-free deformations.
This paper began as an effort to find similar constructions
for other Grassmannians, whose geometry is considerably more complicated 
than that of ${\bf G}_2V$.

This proof of Pieri's formula is an initial step towards
understanding the structure of rational equivalence on these
Grassmann varieties in terms of the combinatorics of the Bruhat
order of the Schubert cellular decomposition. 
A chain in the Bruhat order is a standard skew 
tableau~\cite{Fulton_tableaux}.
Thus the Littlewood-Richardson rule for multiplying two Schubert
classes has an interpretation as a sum over certain chains
in the Bruhat order.
A (as yet unknown) geometric proof of the Littlewood-Richardson
rule for Grassmannians should provide an explanation for this,
similar to what we give for Pieri's formula.

In fact, we believe that all Schubert-type product formulas for
any Grassmannian or flag variety $X$ of any reductive group will
eventually be understood in terms of related combinatorics on
the Bruhat order on Schubert subvarieties on $X$,
perhaps with additional data giving rise to multiplicities.
Some of this picture is already known:
Both Chevalley's~\cite{Chevalley91} formula for multiplication
by hypersurface
Schubert classes and the Pieri-type formulas of
Boe-Hiller~\cite{Hiller_Boe} and 
Pragacz-Ratajski~\cite{Pragacz_Ratajski_Pieri_Odd_I} for
Lagrangian and orthogonal Grassmannians are 
similar to teh form of Pieri's formula (\ref{eq:pieri}), but each has
multiplicities depending upon root system data.
Formulas for multiplying arbitrary Schubert classes in maximal
Lagrangian 
or orthogonal Grassmannians (\cite{Pragacz_S-Q},\cite{Stembridge_shifted})
are similar to the 
Littlewood-Richardson formula, using combinatorics of the
lattice of shifted Young diagrams, the Bruhat order of these
varieties. 
Known formulas for products in the ordinary flag manifold may also
interpreted in terms of chains in the Bruhat order
(\cite{sottile_pieri_schubert},\cite{bergeron_sottile_symmetry}).

It is only in  characteristic zero that general subvarieties of a
Grassmannian intersect generically transversally.
Kleiman~\cite{Kleiman} proves this in characteristic zero and gives a
counterexample in positive characteristic.
In \S\ref{sec:geometry_of_intersections}, we work over an
arbitrary field and give a precise
determination (Theorem~\ref{thm:geometric_intersection}) of when
a special Schubert variety meets a fixed Schubert variety
generically transversally, and describe the components of such an
intersection.  
The geometry of these components is interesting: 
while not an intersection of Schubert varieties, each component is 
`birationally fibred' over such an intersection, with Schubert
variety fibres.
Such cycles are the key to our proof of Pieri's formula
in \S\ref{sec:explicit_rational_equivalences}; 
they are the components of the intermediate cycles in the
deformations used to establish Pieri's formula.

\section{Geometry of Pieri-type 
intersections}\label{sec:geometry_of_intersections}

\subsection{Preliminaries}
Let $k$ be a fixed, but arbitrary, field and $m\leq n$ positive
integers.
Let $V\simeq k^n$ be an $n$-dimensional vector space over $k$
and ${\bf G}_mV$ be the Grassmannian of $m$-planes in $V$.
A {\em complete flag} $\Fdot$ in $V$ is a sequence of subspaces
$$
0=F_{n+1} \subset F_n\subset \cdots\subset F_2\subset F_1 = V
$$
of $V$ where $\dim F_j = n+1-j$.
Let $\Span{S}$ denote the linear span of a subset $S$ of $V$.
We let ${[n]\choose m}$ be the set of all $m$-element subsets of
$[n]:=\{1,2,\ldots,n\}$, considered as decreasing sequences
$\alpha$ of length $m$:
$n\geq\alpha_1>\alpha_2>\cdots>\alpha_m\geq 1$.
A complete flag $\Fdot$ and a sequence $\alpha\in{[n]\choose m}$
together determine a Schubert (sub)variety  of ${\bf G}_mV$,
$$
\Omega_\alpha\Fdot \ :=\ \{H\in{\bf G}_mV\,|\,
\dim H\cap F_{\alpha_j}\geq j, \ 1\leq j\leq m\}.
$$
This variety has codimension $|\alpha|:= \sum \alpha_i-i$.
A {\em special Schubert variety} is the subvariety of all
$m$-planes $H$ which have a nontrivial intersection with
a single subspace $F_{m+s}$ in the flag, $\Omega_{m+s,m-1,\ldots,2,1}\Fdot$.
We use a compact notation for special Schubert varieties.
Let $L:= F_{m+s}$, a subspace of dimension $n+1-m-s$, and 
define 
$$
\Omega_L \ :=\ \Omega_{m+s,m-1,\ldots,2,1}\Fdot.
$$

Two subvarieties meet {\em generically transversally} if they
intersect  trans\-versally along a dense subset of every
component of their  intersection.
They meet {\em improperly} if the codimension of their
(non-empty) intersection is less than the sum of their
codimensions. 
A subspace $L$ meets a flag $\Fdot$ {\em properly} if it meets
each subspace $F_i$ properly.

To simplify some assertions and formulae, we adopt the convention
that if $\gamma$ is a decreasing sequence of length $m$ with
$\gamma_1>n$, then  $\Omega_\gamma\Fdot = \emptyset$.
Similarly, if  the dimension of a subspace is asserted to be 
negative, we intend that subspace to be $\{0\}$.
Also, $\dim \{0\} = -\infty$.

Let $\alpha\in{[n]\choose m}$ and $r$ be a positive integer.
Define $\alpha*r\subset{[n]\choose m}$ to be the set of those
$\beta\in{[n]\choose m}$ with  
$\beta_1\geq\alpha_1>\beta_2\geq\cdots>\beta_m\geq\alpha_m$
and $|\beta|= |\alpha|+r$.
If $\beta\in \alpha*r$, define  
$j(\alpha,\beta)$ to be the first index $i$ where $\beta_i$
differs from $\alpha_i$, $\min\{i\,|\,\beta_i>\alpha_i\}$. 
For $1\leq j\leq m$, let $\delta^j$ be the Kroenecker delta, the
sequence  with a 1 in the $j$th position and $0$'s elsewhere.

\subsection{The cycle $X_\beta(j,F\!\!_{\mbox{\bf .}},L)$}
\label{sec:intermediate_cycle} 
Central to the geometry of Pieri-type
intersections are the components, $X_\beta(j,\Fdot,L)$, of
reducible intersections.
These subvarieties are also components of cycles
intermediate in deformations establishing Pieri's formula.
Let $\beta\in {[n]\choose m}$, $1\leq j\leq m$ be an integer, $\Fdot$ a flag,
and $L$ a linear subspace in $V$.
Define
$$
X_\beta(j,\Fdot,L) \ :=\  \{H\in \Omega_\beta\Fdot\,|\,
\dim H\cap F_{\beta_j}\cap L \geq 1\}, 
$$ 
a subvariety of $\Omega_\beta\Fdot\bigcap \Omega_L$.

The following theorem gives precise conditions on $L$ and
$\Fdot$ which determine whether 
$\Omega_\alpha\Fdot\bigcap \Omega_L$ is improper, generically
transverse, or irreducible.
Moreover, it computes the components of the intersection in the
crucial case of a generically transverse intersection
with the maximal number of irreducible components.

\subsection{Theorem}\label{thm:geometric_intersection}
{\em 
Let $\alpha\in {[n]\choose m}$, $s>0$, $\Fdot$ be a complete
flag in $V$, and  $L\in  {\bf G}_{n+1-m-s}V$.
\begin{enumerate}
\item[(1)] Suppose $\dim F_{\alpha_j}\cap L > n+2-\alpha_j-j-s$
and $F_{\alpha_j}\cap L \neq \{0\}$, for some $1\leq j\leq m$.
Then 
$\Omega_\alpha\Fdot \bigcap \Omega_L$ is improper.
Otherwise, it is generically transverse.
\item[(2)] Suppose $\dim  F_{\alpha_j}\cap L = n+2-\alpha_j-j-s$ for 
each $1\leq j\leq m$.  Let $\Mdot$ be any flag satisfying
$M_{\alpha_j} = F_{\alpha_j}$  and 
$M_{\alpha_j+1} \supset \Span{F_{\alpha_{j-1}}, F_{\alpha_j}\cap L}$,
for $1\leq j\leq m$.
Then  $\Omega_\alpha \Fdot$ meets $\Omega_L$ generically transversally, and 
$$
\Omega_\alpha\Fdot \bigcap \Omega_L\ =\ 
\sum_{\beta\in \alpha*1} X_\beta(j(\alpha,\beta),\Mdot,L).
$$
\item [(3)] Suppose $\dim F_{\alpha_j}\cap L< n+2-\alpha_j-j-s$
for  each $1\leq j<m$ and 
$F_{\alpha_m}$ meets $L$ properly, so that 
$\dim F_{\alpha_m}\cap L = n+2-\alpha_m-m-s$.
Then $\Omega_\alpha\Fdot \bigcap \Omega_L$ is irreducible.
\end{enumerate}
}
\smallskip

Note that $n+2-\alpha_j-j-s$, the critical dimension for
$F_{\alpha_j}\cap L$ in this theorem, exceeds the expected dimension of
$n+2-\alpha_j-m-s$ by $m-j$.
Thus, it is not necessary for $\Fdot$ and $L$ to meet properly for
$\Omega_\alpha\Fdot\bigcap\Omega_L$ to be generically transverse or
even irreducible.
However, it is necessary that $F_{\alpha_m}$ and $L$ meet properly.
We also see that, as the relative position of $\Fdot$ and $L$ becomes more
degenerate, the intersection $\Omega_\alpha\Fdot\bigcap\Omega_L$
`branches' into components, one for each
$j$ such that $\dim F_{\alpha_j}\cap L = n+2-\alpha_j-j-s$, and it
will attain excess intersection if 
$\dim F_{\alpha_j}\cap L > n+2-\alpha_j-j-s$, for even one $j$.

\subsection{Remark}\label{remarkI}
In the situation of Theorem~\ref{thm:geometric_intersection}(2), 
if $\beta\in \alpha*1$ and $j(\alpha,\beta) =1$,
then $\beta = \alpha+\delta^1$.
Suppose further that  $M_{\alpha_1}\cap L= M_{\alpha_1+s}$.
Then 
$$
X_\beta(1,\Mdot,L) \ 
=\  \Omega_{\alpha+s\delta^1}\Mdot\ =\ \Omega_{\beta+(s-1)\delta^1}\Mdot,
$$
so we have 
$$
\Omega_\alpha\Fdot \bigcap \Omega_L\ =\ 
\sum_{\stackrel{\mbox{\scriptsize${\beta}{\in}{\alpha}{*}{1}$}}
{j(\alpha,\beta)=1}} \Omega_{\beta+(s-1)\delta^1
}\Mdot\ 
+\sum_{\stackrel{\mbox{\scriptsize${\beta}{\in}{\alpha}{*}{1}$}}
{j(\alpha,\beta)>1}} 
X_\beta(j(\alpha,\beta),\Mdot,L).
$$

We prove Theorem~\ref{thm:geometric_intersection} in
\S\ref{sec:proof_geometric_intersection}.
First, we study the varieties $X_\beta(j,\Fdot,L)$.
Let $\beta\in{[n]\choose m}$, $\Fdot$ be a complete flag, and
$1\leq j\leq m$ an integer.
The map from $\Omega_\beta\Fdot$ to ${\bf G}_jF_{\beta_j}$ given
by $H\longmapsto H\cap F_{\beta_j}$ is only defined on the dense
locus in $\Omega_\beta\Fdot$ of those $H$ where  $\dim H\cap
F_{\beta_j}=j$. 
Resolving the ambiguity of this map gives the variety
$$
\widetilde{\Omega_\beta}^j\!\Fdot\ :=\ 
\{(H,K) \in \Omega_\beta\Fdot \times {\bf G}_jF_{\beta_j}
\,|\, K\subset H\mbox{\ \ and\ \ } 
\dim K\cap F_{\beta_i}\geq i,
\ 1\leq i\leq j\}.
$$
In Lemma~\ref{lemma:rational_fibration}, we show that the 
projection to ${\bf G}_jF_{\beta_j}$
realizes $\widetilde{\Omega_\beta}^j\!\Fdot$ as a fibre bundle
with base and fibres themselves Schubert varieties.
The following definitions are needed to describe the base and
fibres. 
Let $p$ be the first projection and $\pi$ the second.
For $K\subset V$, let $\Fdot/K$ be the image of the flag $\Fdot$
in $V/K$. 
Let $\Fdot|_{\beta_j}$ be the flag 
$$
F_{\beta_j}\supset F_{\beta_j +1}\supset\cdots\supset F_n
$$
and  $\beta|_j\in {[n+1-\beta_j]\choose j}$ the sequence
$$
\beta_1-\beta_j+1>\cdots>\beta_{j-1}-\beta_j+1>1=
(\beta|_j)_j.
$$  
Unraveling this definition shows  
$(\Fdot|_{\beta_j})_{(\beta|_j)_i}=F_{\beta_i}$, for $i\leq j$. 

\subsection{Lemma.}\label{lemma:rational_fibration} {\em
Let $\beta\in {[n]\choose m}$, $\Fdot$ be a flag, and 
$1\leq j\leq m$.
Then $p$ is an isomorphism over the dense subset
$\{ H\in \Omega_\beta\Fdot\, |\, \dim H\cap F_{\beta_j} = j\}$.
Also, $\pi$ exhibits 
$\widetilde{\Omega_\beta}^j\!\Fdot$ as a fibre bundle with 
base $\Omega_{\beta|_j}\Fdot|_{\beta_j}$ whose fibre over 
$K\in\Omega_{\beta|_j}\Fdot|_{\beta_j}$ is 
the Schubert variety 
$\Omega_{\beta_{j+1}\ldots\,\beta_m}\Fdot/K\subset{\bf G}_{m-j}V/K$.
Moreover, each fibre of $\pi$ meets the locus where $p$ 
is an isomorphism.
}\smallskip

\noindent{\bf Proof:}
We describe the fibres of $\pi$.
Note that Schubert varieties 
have a dual description:
$$
H\in \Omega_\beta\Fdot \ \Longleftrightarrow\ 
\dim\frac{H}{H\cap F_{\beta_i}} \leq m-i,\ \mbox{\ for\ } 1\leq i\leq m.
$$

If $K\in \Omega_{\beta|_j}\Fdot|_{\beta_j}$, 
then $K\subset F_{\beta_j}\subset F_{\beta_i}$, for $i>j$.
Thus $\left( \Fdot/K\right)_{\beta_i} = F_{\beta_i}/K$, for $i>j$.
Hence, if $H$ is in the fibre over $K$, then $H\in \Omega_\beta\Fdot$ and 
$K\subset H$, so
$$
\dim \frac{H/K}{H/K\cap \left( \Fdot/K\right)_{\beta_i}} 
\ =\ \dim\frac{H}{H\cap F_{\beta_i}} \leq m-i,\  \mbox{\ for\ }
j< i\leq m.
$$
Thus  $H/K \in \Omega_{\beta_{j+1}\ldots\,\beta_m}\Fdot/K$.
The reverse implication is similar
and the remaining assertions follow easily from the definitions.
\QED

Reformulating the definition of $X_\beta(j,\Fdot,L)$ 
in these terms gives a useful characterization.

\subsection{Corollary.}\label{cor:intermediate_cycles}
$
X_\beta(j,\Fdot,L) \ =\ p(\pi^{-1}( 
\Omega_{\beta|_j}\Fdot|_{\beta_j} \bigcap \Omega_{F_{\beta_j}\cap L})).
$
\smallskip

Since the fibres of $\pi$ meet the locus where $p$ is an isomorphism,
the map 
$$
p \ :\ \pi^{-1}( 
\Omega_{\beta|_j}\Fdot|_{\beta_j} \bigcap \Omega_{F_{\beta_j}\cap L})
\longrightarrow X_\beta(j,\Fdot,L)
$$
is  proper and birational.
Thus, while $X_\beta(j,\Fdot,L)$ is neither a Schubert variety nor an 
intersection of Schubert varieties, it is `birationally fibred' over 
an intersection of Schubert varieties with Schubert variety fibres,
and hence is intermediate between these extremes.

\subsection{Tangent spaces to Schubert varieties}
Let $H\in {\bf G}_m V$ and $K\in {\bf G}_{n-m}V$ be complementary
subspaces, so $H\cap K = \{0\}$.
The open set $U\subset {\bf G}_m V$ of those $H'$ with  $H'\cap K = \{0\}$
is identified with $\Hom(H,K)$ by 
$\phi \in \Hom(H,K) \mapsto \Gamma_\phi$, the graph of $\phi$ in 
$H\oplus K =V$.
This shows the tangent space of ${\bf G}_mV$ at $H$,  $T_H{\bf G}_mV$,
is equal to $\Hom(H,V/H)$, since $K$ is canonically isomorphic to $V/H$.
The intersection of a Schubert variety $\Omega_\alpha\Fdot$ containing $H$
with this open set $U$ can be used to determine whether $\Omega_\alpha\Fdot$
is smooth at $H$ and its tangent space at $H$.
This gives the following description:
If $H\in {\bf G}_mV$ and  $\dim H\cap F_{\alpha_j} = j$ for 
$1\leq j\leq m$, then $\Omega_\alpha\Fdot$ is smooth at $H$ and
$$
T_H\Omega_\alpha\Fdot\  =\ \{ \phi\in \Hom(H,V/H)\,|\,
\phi(H\cap F_{\alpha_j}) \subset (F_{\alpha_j} +H)/H,\ 1\leq j\leq m\}.
$$
Similarly, if  $H\in {\bf G}_mV$, 
$L\in {\bf G}_{n+1-m-s}V$, and $\dim H\cap L = 1$,
then   $\Omega_L$ is smooth at $H$ and the tangent space of $\Omega_L$
at $H$ is 
$$
T_H\Omega_L \ =\ \{\phi\in \Hom(H,V/H)\,|\,
\phi(H\cap L)\subset (L+H)/H\}.
$$

Let $P$ be the subgroup of $GL(V)$ stabilizing the partial flag
$F_{\alpha_1}\subset F_{\alpha_2}\subset \cdots\subset F_{\alpha_m}$.
The orbit  $P\cdot L'$ consists of those $L$ with 
$\dim F_{\alpha_j}\cap L=\dim F_{\alpha_j}\cap L'$ for 
$1\leq j\leq m$.
Similarly,  $L\in \overline{P\cdot L'}$
if $\dim F_{\alpha_j}\cap L\geq \dim F_{\alpha_j}\cap L'$ for 
$1\leq j\leq m$.
If
$P\cdot L= P\cdot L'$, then 
$\Omega_\alpha\Fdot\bigcap\Omega_L 
\simeq\Omega_\alpha\Fdot\bigcap\Omega_{L'}$.
Thus $P$-orbits on ${\bf G}_{n+1-m-s}V$ determine the isomorphism type
of Pieri-type intersections.

\subsection{Lemma} \label{lemma:P-orbits}
{\em
Suppose that $L,L'\in{\bf G}_{n+1-m-s}V$  with 
$L\in \overline{P\cdot L'}$.
Then
\begin{enumerate}
\item [(1)] $\dim\Omega_\alpha\Fdot\bigcap\Omega_L \geq 
\dim\Omega_\alpha\Fdot\bigcap\Omega_{L'}$.
\item [(2)] If $\,\Omega_\alpha\Fdot\bigcap\Omega_L$ is 
generically transverse, then 
 $\Omega_\alpha\Fdot\bigcap\Omega_{L'}$ is generically transverse.
\item [(3)] If $\,\Omega_\alpha\Fdot\bigcap\Omega_L$ is generically 
transverse and 
irreducible, then  $\Omega_\alpha\Fdot\bigcap\Omega_{L'}$ is generically 
transverse and irreducible.
\end{enumerate} 
}\smallskip

\noindent{\bf Proof:}
Let $\psi:{\bf P}^1 \rightarrow \overline{P\cdot L'}$ be a map with 
$\psi(0)=L$ and $\psi({\bf P}^1) \cap (P\cdot L') \neq \emptyset$.
Then $\Omega_\alpha\Fdot\bigcap \Omega_{\psi(t)}$ is isomorphic to 
$\Omega_\alpha\Fdot\bigcap\Omega_{L'}$, for any
$t\in \psi^{-1}(P\cdot L')$.
The lemma  follows by considering the subvariety 
of ${\bf P}^1\times {\bf G}_mV$ whose fibre over $t\in {\bf P}^1$ 
is $\Omega_\alpha\Fdot\bigcap \Omega_{\psi(t)}$.
\QED

\subsection{\bf Proof of
Theorem~\ref{thm:geometric_intersection}:}
\label{sec:proof_geometric_intersection}
Let $\alpha\in{[n]\choose m}$, $s>0$, $\Fdot$ be a complete flag, and
$L\in {\bf G}_{n+1-m-s}V$. 
The conditions on $L$ in statement (2), that 
$\dim F_{\alpha_j}\cap L = n+2-\alpha_j-j-s$ for each $j$, determine a
$P$-orbit, 
which is the closure of any $P$-orbit $P\cdot L'$, where 
$\dim F_{\alpha_j}\cap L' \leq n+2-\alpha_j-j-s$ for each $j$.
Thus (2) and Lemma~\ref{lemma:P-orbits}(2) together imply that if
$\dim F_{\alpha_j}\cap L \leq n+2-\alpha_j-j-s$ for each $j$, then
$\Omega_\alpha\Fdot\bigcap\Omega_L$ is 
generically transverse, proving the second part of (1).

For the first part of (1), suppose 
$\dim F_{\alpha_j}\cap L > n+2-\alpha_j-j-s$ and let 
$L':= F_{\alpha_j}\cap L\neq \{0\}$.
Then $L'$ has codimension at most $j+s-1$ in $F_{\alpha_j}$.
Hence 
$\Omega_{\alpha|_j}\Fdot|_{\alpha_j} \bigcap \Omega_{L'}\neq \emptyset$
and so has  codimension in $\Omega_{\alpha|_j}\Fdot|_{\alpha_j}$ at most that 
of $\Omega_{L'}$ in ${\bf G}_jF_{\alpha_j}$, which is at most $s-1$.
Thus 
$$
X_\alpha(j,\Fdot,L)\ = \ 
p (\pi^{-1}(\Omega_{\alpha|_j}\Fdot|_{\alpha_j} \bigcap \Omega_{L'}))
$$
which has codimension less than $s$ in $\Omega_\alpha\Fdot=
p (\pi^{-1}(\Omega_{\alpha|_j}\Fdot|_{\alpha_j}))$.
Hence  $\Omega_\alpha\Fdot\bigcap \Omega_L$ 
is improper, as 
$X_\alpha(j,\Fdot,L)\subset \Omega_\alpha\Fdot\bigcap \Omega_L$, proving (1).
\smallskip

Before proceeding with the rest of the proof, we make a computation.
Suppose $\dim F_{\alpha_j}\cap L \leq n+2-\alpha_j-j-s$ for 
$1\leq j\leq m$ and $F_{\alpha_m}\cap L \not\subset F_{\alpha_{m-1}}$.
Then there exists $H\in \Omega_\alpha\Fdot\bigcap \Omega_L$ with 
$\dim H\cap F_{\alpha_j}=j$ for $1\leq j\leq m$, $\dim H\cap L = 1$,
and $H\cap L\not\subset F_{\alpha_{m-1}}$:
Inductively choose linearly independent vectors $f_j\in F_{\alpha_j}$
for $1\leq j\leq m$ as follows.
Let $f_1\in F_{\alpha_1} - \{0\}$.
Then for $1<j<m$ suppose that $f_1,\ldots,f_{j-1}$ have been chosen.
Since $\dim F_{\alpha_j}$ exceeds
$$
\dim F_{\alpha_j}\cap \Span{L,f_1,\ldots,f_{j-1}}
\ \leq\  n+2-\alpha_j-j-s+(j-1)\  =\ n+1-\alpha_j-s
$$
and $F_{\alpha_{j-1}}\not\subset F_{\alpha_j}$, we can select a vector
$f_j$ in 
$$
 F_{\alpha_j}-\Span{L,f_1,\ldots,f_{j-1}}- F_{\alpha_{j-1}}.
$$
Let $f_m\in F_m\cap L-F_{\alpha_{m-1}}$, and let 
$H := \Span{f_1,\ldots,f_m}$.
Then $H\in \Omega_\alpha\Fdot\bigcap\Omega_L$, $\dim H\cap F_{\alpha_j}=j$ for
$1\leq j\leq m$, $\dim H\cap L=1$, and $H\cap L\not\subset F_{\alpha_{m-1}}$.
Let $X_m^\circ$ be the set of all such $H$.
For $H\in X_m^\circ$,
$$
T_H\Omega_\alpha\Fdot\bigcap T_H\Omega_L\ =\ \{ \phi\in T_H\Omega_\alpha\Fdot
\,|\, \phi(H\cap L) \subset (F_{\alpha_m}\cap L +H)/H\}.
$$
This has codimension in $T_H\Omega_\alpha\Fdot$ equal to
$\dim(F_{\alpha_m}+H) - 
\dim(F_{\alpha_m}\cap L + H) =  s$.
Thus $\Omega_\alpha\Fdot$ and $\Omega_L$ meet transversally along 
$X_m^\circ$.
\smallskip

We show  (2).
Suppose $\dim F_{\alpha_j}\cap L = n+2-\alpha_j-s$ for each 
$1\leq j\leq m$.
Let $\Mdot$ be any flag satisfying $M_{\alpha_j}=F_{\alpha_j}$ and 
$M_{\alpha_j +1}\supset\Span{F_{\alpha_{j-1}},\,F_{\alpha_j}\cap L}$,
for $1\leq j\leq m$.
Let $H\in \Omega_\alpha\Fdot\bigcap \Omega_L$.
Then there is some $1\leq j\leq m$ with 
$H\cap L\cap F_{\alpha_j}\not\subset F_{\alpha_{j-1}}$.
Since $\dim H\cap F_{\alpha_{j-1}}\geq j-1$, we have
$\dim H\cap\Span{F_{\alpha_{j-1}}, F_{\alpha_j}\cap L}\geq j$
and so $\dim H\cap M_{\alpha_j+1}\geq j$.
Thus $H\in \Omega_{\alpha+\delta^j}\Mdot$, if
$\alpha+\delta^j\in{[n]\choose m}$.
But this is the case, as $\alpha_j+1<\alpha_{j-1}$, for otherwise 
dimensional considerations imply that 
$L\cap F_{\alpha_j}=L\cap F_{\alpha_{j-1}}\subset F_{\alpha_{j-1}}$.
Let $\beta := \alpha+\delta^j \in \alpha*1$.
Then $j(\alpha,\beta) = j$
and $H\in X_\beta(j(\alpha,\beta),\Mdot,L)$, since 
$H\in \Omega_\beta \Mdot$ and $\dim H\cap L\cap M_{\beta_j}\geq 1$.
Conversely, if $\beta\in \alpha*1$,  then 
$\Omega_\beta\Mdot\subset\Omega_\alpha\Fdot$, so 
$X_\beta(j(\alpha,\beta),\Mdot,L)\subset\Omega_\alpha\Fdot\bigcap\Omega_L$.  
This shows 
$$
 \Omega_\alpha\Fdot\bigcap \Omega_L\ =\ \sum_{\beta\in \alpha*1}
X_\beta(j(\alpha,\beta),\Mdot,L).
$$
We claim this intersection is generically transverse.
Let $\beta\in \alpha*1$ and $j := j(\alpha,\beta)$. 
Then 
$X_\beta(j,\Mdot,L)$ has an open subset $X^\circ_j$ 
consisting of those $H$ with  $\dim H\cap F_{\alpha_i} = i$ for 
$1\leq i\leq m$, $\dim H\cap L= 1$, and $H\cap L\subset F_{\alpha_j}$
but $H\cap L\not\subset F_{\alpha_j-1}$.
As with $X^\circ_m$ above, $X^\circ_j$ is nonempty, so it is a dense open
subset of $X_\beta(j,\Mdot,L)$
For $H\in X^\circ_j$,
$$
T_H\Omega_\alpha\Fdot \bigcap T_H\Omega_L \ = \ 
\{\phi\in T_H\Omega_\alpha\Fdot \,|\,
\phi(H\cap L) \subset (L\cap F_{\alpha_j} +H)/H\}.
$$
Since 
$\dim(F_{\alpha_j} +H)- \dim(L\cap F_{\alpha_j} +H)= s$,
this has  codimension $s$ 
in $T_H\Omega_\alpha\Fdot$, showing that $\Omega_\alpha\Fdot$ and 
$\Omega_L$ meet transversally along $X^\circ_j$, a dense 
subset of $X_\beta(j(\alpha,\beta),\Mdot,L)$.%
\smallskip

By Lemma~\ref{lemma:P-orbits}(3),  it suffices to prove 
a special case of (3):
\begin{enumerate}
\item[(3)$'$] {\em 
If  $F_{\alpha_m}$ meets $L$ properly, and 
for $1\leq j<m$, $\dim F_{\alpha_j}\cap L= n+2-\alpha_j-j-(s+1)$, 
then $\Omega_\alpha\Fdot \bigcap \Omega_L$ is irreducible.}
\end{enumerate}
These conditions imply  $F_{\alpha_m}\cap L \not\subset F_{\alpha_{m-1}}$.
In the notation of \S \ref{remarkI}, let $L':=F_{\alpha_{m-1}}\cap L$, 
$\Fpdot := \Fdot|_{\alpha_{m-1}}$,   and $\alpha':= \alpha|_{m-1}$.
Consider
$$
X_\alpha(m-1,\Fdot,L) \ =\  p(\pi^{-1}(
\Omega_{\alpha'}\Fpdot\bigcap 
\Omega_{L'})).
$$
For $j\leq m-1$, 
\begin{eqnarray*}
\dim F_{\alpha_j}\cap L' &=& n+2-\alpha_j-j-(s+1)\\
&=& \dim F_{\alpha_{m-1}} +2 -\alpha'_j -j-(s+1),
\end{eqnarray*}
so 
$L'$ and $\Fpdot$ satisfy the conditions of (2) for the pair 
$\alpha', s+1$.
Thus $\Omega_{\alpha'}\Fpdot\bigcap \Omega_{L'}$
is generically transverse, which implies that 
$X_\alpha(m-1,\Fdot,L)$ has codimension $s+1$ in 
$\Omega_\alpha\Fdot$ and hence is a proper subvariety of 
$\Omega_\alpha\Fdot\bigcap\Omega_L-X_\alpha(m-1,\Fdot,L)$.
Since  $X^\circ_m$ is dense in 
$\Omega_\alpha\Fdot\bigcap\Omega_L-X_\alpha(m-1,\Fdot,L)$, this  completes
the proof of (3)$'$.
\QED

\section{Construction of explicit rational 
equivalences}\label{sec:explicit_rational_equivalences}

Theorem~\ref{thm:geometric_intersection} shows that for $L$ in 
a dense subset of ${\bf G}_{n+1-m-b}V$, the intersection 
$\Omega_\alpha\Fdot\bigcap \Omega_L$ is generically transverse
and irreducible.
We use Theorem~\ref{thm:geometric_intersection}(2)
to study  such a cycle as $L$ `moves out of'
this set, ultimately deforming it into the cycle
$\sum_{\gamma\in \alpha*b}\Omega_\gamma\Fdot$.

\subsection{Families and Chow varieties}\label{sec:Chow}
Suppose $\Sigma \subset ({\bf P}^1-\{0\}) \times {\bf G}_mV$ 
has equidimensional
fibres over ${\bf P}^1-\{0\}$.
Then its Zariski closure $\overline{\Sigma}$ 
in ${\bf P}^1\times{\bf G}_mV$ has equidimensional fibres over
${\bf P}^1$.
Denote the fibre of $\overline{\Sigma}$ over 
$0$ by $ \lim_{t\rightarrow 0} \Sigma_t$,
where $\Sigma_t$ is the fibre of $\Sigma$ over $t\in {\bf P}^1-\{0\}$.
The association of a point $t$ of ${\bf P}^1$ 
to the fundamental cycle of the 
fibre $\overline{\Sigma}_t$ determines a morphism 
${\bf P}^1\rightarrow \Chow {\bf G}_mV$.
Moreover, if $\Sigma$ is  defined over $k$, then so 
is the map ${\bf P}^1\rightarrow \Chow {\bf G}_mV$ 
(\cite{Samuel}, \S I.9).

\subsection{The cycle $Y_{\alpha,r}(F\!\!_{\mbox{\bf .}},L)$}
\label{sec:Y_alpha,r}
In \S \ref{sec:intermediate_cycle}, we defined the components
$X_\beta(\alpha,\Fdot,L)$ of the cycles intermediate between 
$\Omega_\alpha\Fdot\bigcap \Omega_L$ and 
$\sum_{\gamma\in \alpha*b}\Omega_\gamma\Fdot$.
Here, we define those intermediate cycles, 
$Y_{\alpha,r}(\Fdot,L)$, which 
are parameterized by subspaces $L$ in certain Schubert
cells $U_{\alpha,s}\Fdot$ of ${\bf G}_{n+1-m-s}$.
Let $U_{\alpha,s}\Fdot$ be the set  of those
$L\in G_{n+1-m-s}V$ such that
\begin{enumerate}
\item[(1)] $F_{\alpha_1} \cap L = F_{\alpha_1+s}$, and 
\item[(2)] $F_{\alpha_j}\cap L = F_{\alpha_j+1}\cap L$, and has dimension 
$n+2-\alpha_j-j-s$, for $1\leq j\leq m$.
\end{enumerate}
These conditions are consistent and determine 
$\dim F_i\cap L$ for $1\leq i\leq n$.
For example,
$$
\mbox{(\ref{sec:Y_alpha,r})}\hspace{.57in} 
\alpha_j<i<\alpha_{j-1}\ \ \Longrightarrow\ \ 
\dim F_i\cap L \ =\ \dim F_i + 1-j-(s-1).
\hspace{.87in}
$$
Thus $U_{\alpha,s}\Fdot$ is a single Schubert cell of $G_{n+1-m-s}V$.
Specifically, $U_{\alpha,s}\Fdot$
is the dense cell of $\Omega_\beta\Fdot$,
where $\beta\in{[n]\choose n+1-m-s}$ is defined as follows:
If $\alpha_1\leq n+1-s$, then 
$\beta=[n]-\alpha-\{\alpha_1+1,\ldots,\alpha_1+s-1\}$.
Otherwise, $\beta$ is the smallest $n+1-m-s$ integers in 
$[n]-\alpha$.

For $\beta\in \alpha*r$, recall that 
$j(\alpha,\beta) = \min\{i\,|\, \alpha_i<\beta_i\}$.
If $L\in U_{\alpha,s}\Fdot$, define the cycle
$$
Y_{\alpha,r}(\Fdot,L)\  := 
\sum_{\stackrel{\mbox{ \scriptsize ${\beta}{\in}{\alpha}{*}{r}$}}
{j(\alpha,\beta)=1}}
\Omega_{\beta+ (s-1)\delta^1}\Fdot
\ +
\sum_{\stackrel{\mbox{ \scriptsize ${\beta}{\in}{\alpha}{*}{r}$}}
{j(\alpha,\beta)>1}}
X_\beta(j(\alpha,\beta),\Fdot,L).
$$
Let ${\cal G}_{\alpha,s,r}\Fdot \subset \Chow {\bf G}_mV$ be the set of 
these cycles $Y_{\alpha,r}(\Fdot,L)$ for 
$L\in U_{\alpha,s}\Fdot$.
Since $U_{\alpha,s}\Fdot$ is a Schubert cell,
${\cal G}_{\alpha,s,r}\Fdot$ 
is an orbit of the Borel subgroup 
stabilizing $\Fdot$ and hence is rational.

\subsection{Remark}\label{remarkII}
Suppose  $L\in U_{\alpha,s}\Fdot$, then by 
 Remark~\ref{remarkI}, 
\begin{eqnarray*}
\Omega_\alpha\Fdot\bigcap \Omega_L &=&
\sum_{\stackrel{\mbox{ \scriptsize ${\beta}{\in}{\alpha}{*}{1}$}}
{j(\alpha,\beta)=1}}
\Omega_{\beta+ (s-1)\delta^1}\Fdot\ 
+
\sum_{\stackrel{\mbox{ \scriptsize ${\beta}{\in}{\alpha}{*}{1}$}}
{j(\alpha,\beta)>1}}
X_\beta(j(\alpha,\beta),\Fdot,L)\\
&=& Y_{\alpha,1}(\Fdot,L).
\end{eqnarray*}

The following lemma parameterizes our explicit rational 
equivalences.
It is identical to   Lemma~6.1 of~\cite{sottile_real_lines}.

\subsection{Lemma.}\label{lemma:limits_are_good}
{\em 
Let $l\leq n$ and let 
$\Mdot$ be a complete flag in $M\simeq k^n$.
Suppose $L_{\infty}$ is a hyperplane containing $M_l$ but not $M_{l-1}$.
Then there exists a pencil of hyperplanes $L_t$, for 
$t\in {\bf P}^1$,
such that if $t\neq 0$, then $L_t$ contains  $M_l$ but not $M_{l-1}$
and,  for each $i\leq l-1$, the family of codimension $i+1$ planes
induced by $M_i\cap L_t$ for $t\neq 0$ has fibre $M_{i+1}$ over $0$.
}
\smallskip

\noindent{\bf Proof:}
Let $x_1,\ldots,x_n$ be a basis of $M^*$ such that 
$L_{\infty}= \Span{x_{l-1}}^{\perp}$ and
$M_i = \Span{x_1,\ldots,x_{i-1}}^\perp$.
Let $e_1,\ldots,e_n$ be a basis for $M$ dual to $x_1,\ldots,x_n$
and define 
$$
L_t \ :=\  \Span{M_l, te_j + e_{j+1}\,|\, 1\leq j\leq l-2}.
$$
For $t\neq 0$  and $1\leq i\leq l-1$, 
$M_i \cap L_t = \Span{M_l, te_j + e_{j+1}\,|\, i\leq j\leq l-2}$
and so has dimension $n-i$.
The fibre of this family at $t=0$ is 
$\Span{M_l, e_{j+1}\,|\, i\leq j\leq l-2} = M_{i+1}$.
\QED

\subsection{Theorem.} \label{thm:inductive_engine}{\em
Let $\alpha\in {[n]\choose m}$, s,r be positive integers and $\Fdot$ a 
flag in $V$.
Let $M\in U_{\alpha,s-1}\Fdot$ and define $\Mdot$ to be the flag in 
$M$ consisting of the subspaces
in $\Fdot\cap M$.

Let $L_\infty\subset M$ be any hyperplane containing $F_{\alpha_1+s}$
but not $F_{\alpha_1+s-1}$.
Suppose $L_t$ is the family of hyperplanes of $M$
given by Lemma~\ref{lemma:limits_are_good}.
Then
\begin{enumerate}
\item[(1)] For $t\neq 0$, $L_t \in U_{\alpha,s}\Fdot$.
\item[(2)] ${\displaystyle \lim_{t\rightarrow 0} Y_{\alpha,r}(\Fdot,L_t)
= Y_{\alpha, r+1}(\Fdot,M)}$.
\end{enumerate}
}

\subsection{Theorem.}\label{thm:pieri} {\bf [Pieri's Formula]} 
{\em Let $\alpha\in{[n]\choose m}$, $\Fdot$ be a complete flag in $V$, and
$K\in{\bf G}_{n+1-m-b}V$ be a subspace which meets $\Fdot$ properly.
Then the cycle $\Omega_\alpha\Fdot\bigcap\Omega_K$, a generically
transverse intersection, is rationally equivalent to 
$
\sum_{\gamma\in \alpha*b}\Omega_\gamma\Fdot$. 
Thus, in the Chow ring $A^*{\bf G}_mV$ of ${\bf G}_mV$, 
$$
[\Omega_\alpha\Fdot]\cdot[\Omega_K] \ =\ 
\sum_{\gamma\in \alpha*b}[\Omega_\gamma\Fdot].
$$

Moreover, let ${\cal G}\subset \Chow {\bf G}_mV$ be the set of cycles arising 
as generically transverse intersections of the form
$\Omega_\alpha\Fdot\bigcap\Omega_K$ for $K\in{\bf G}_{n+1-m-b}V$.
Then  one may give $b+1$ explicit rational deformations inducing this
rational equivalence, where the cycles at the $i$th stage are of the form 
$Y_{\alpha,i}(\Fdot,M)$, with $M\in U_{\alpha,b+1-i}\Fdot$, and all are
within $\overline{\cal G}$.
}\smallskip

The Borel subgroup $B$ of $GL(V)$ stabilizing $\Fdot$ also 
stabilizes ${\cal G}$ and the cycle 
$\sum_{\gamma\in\alpha*b}\Omega_\gamma\Fdot$
is the only $B$-stable cycle in this component of the Chow variety.
As Hirschowitz~\cite{Hirschowitz} observed, 
this implies $\sum_{\gamma\in\alpha*b}\Omega_\gamma\Fdot$ is in the Zariski
closure of  ${\cal G}$.
Thus Theorem~\ref{thm:pieri} is an improvement in that 
deformations inducing the rational equivalence are given explicitly.

\subsection{Proof of Pieri's formula using Theorem~3.5}\label{sec:pieri_proof}
Let $b>0$, and $\alpha\in {[n]\choose m}$.
For $1\leq i\leq b$, let $U_i := U_{\alpha,b+1-i}\Fdot$
and ${\cal G}_i := {\cal G}_{\alpha,b+1-i,i}\Fdot$.
Let $U_0\subset {\bf G}_{n+1-m-b}V$ be the (dense) set of those $L$
which meet $F_{\alpha_m}$ properly and for $1\leq j<m$, 
$\dim F_{\alpha_j}\cap L< n+2-\alpha_j-j-b$.
By Theorem~\ref{thm:geometric_intersection},
if $L\in {\bf G}_{n+1-m-b}V$, then 
$\Omega_\alpha\Fdot\bigcap \Omega_L$ is generically transverse and 
irreducible if and only if $L\in U_0$.
Let ${\cal G}_0\subset \Chow {\bf G}_mV$ be the set of cycles
$\Omega_\alpha\Fdot\bigcap \Omega_L$ for $L\in U_0$.

Let $L\in U_b$ and consider the cycle 
$Y_{\alpha,b}(\Fdot,L)\in {\cal G}_b$:
$$
Y_{\alpha,b}(\Fdot,L)\  = 
\sum_{\stackrel{\mbox{ \scriptsize ${\beta}{\in}{\alpha}{*}{b}$}}
{j(\alpha,\beta)=1}}
\Omega_{\beta}\Fdot
\ +
\sum_{\stackrel{\mbox{ \scriptsize ${\beta}{\in}{\alpha}{*}{b}$}}
{j(\alpha,\beta)>1}}
X_\beta(j(\alpha,\beta),\Fdot,L).
$$
We claim 
$Y_{\alpha,b}(\Fdot,L)=\sum_{\beta\in\alpha*b}\Omega_\beta\Fdot$,
the cycle $Y_{\alpha,b}\Fdot$ of the Introduction.
It suffices to show $X_\beta(j(\alpha,\beta),\Fdot,L)=\Omega_\beta\Fdot$
for $\beta\in \alpha*b$ with $j(\alpha,\beta)>1$.
Suppose $j = j(\alpha,\beta)>1$, then 
$$
X_\beta(j,\Fdot,L) \ =\  p(\pi^{-1}(
\Omega_{\beta|_j}\Fdot|_{\beta_j}\bigcap \Omega_{F_{\beta_j}\cap L})).
$$
By formula (\ref{sec:Y_alpha,r}), 
$\dim F_{\beta_j}\cap L = \dim F_{\beta_j}-j+1$, 
as $\alpha_j<\beta_j<\alpha_{j-1}$ and $s=1$.
So $\Omega_{F_{\beta_j}\cap L}$ is ${\bf G}_jF_{\beta_j}$,
since any $j$-plane in $F_{\beta_j}$ meets $F_{\beta_j}\cap L$ 
non-trivially.
Thus $X_\beta(j(\alpha,\beta),\Fdot,L)=\Omega_\beta\Fdot$, by the definition
of $p$ and  $\pi$ in~\S\ref{remarkI}.
\smallskip

Let ${\cal G}\subset \Chow {\bf G}_mV$ be the set of all cycles
$\Omega_\alpha\Fdot \cap \Omega_L$, where 
$L\in {\bf G}_{n+1-m-b}V$ and the intersection is generically transverse.
Then by Theorem~\ref{thm:geometric_intersection} and Remark~\ref{remarkII},
both ${\cal G}_0$ and ${\cal G}_1 $ are subsets of ${\cal G}$.
Arguing as in the proof of Lemma~\ref{lemma:P-orbits} shows 
${\cal G}\subset \overline{{\cal G}_0}$.
Theorem~\ref{thm:inductive_engine} implies
${\cal G}_i\subset \overline{{\cal G}_{i-1}}$ for $2\leq i\leq b$, 
so in particular, $Y_{\alpha,b}\Fdot\in{\cal G}_b\subset \overline{\cal G}$.
Since ${\cal G}_0$ and hence $\overline{\cal G}$ is rational,
$Y_{\alpha,b}\Fdot$ is rationally
equivalent to any cycle in   ${\cal G}$, including
$\Omega_\alpha\Fdot\bigcap\Omega_K$, proving Pieri's formula. 
\smallskip

More explicitly,  one may construct a sequence of parameterized rational 
curves
$\phi_i: {\bf P}^1 \rightarrow \overline{\cal G}$ for 
$1\leq i\leq b$ witnessing this
rational equivalence.
For $2\leq i\leq b$, select subspaces $M_i\in U_i$ and pencils 
$L_{i,t}$ of hyperplanes of $M_i$ by downward induction on $i$ as follows:
Choose $M_b\in U_b$.
Given $M_i\in U_i$, let  $L_{i,t}$ 
be a pencil of hyperplanes of $M_i$ as in 
Theorem~\ref{thm:inductive_engine}, let $M_{i-1}: = L_{i,\infty}$, and 
continue.
Then for each $i$, if $t\neq 0$, $L_{i,t}\in U_{i-1}$.
Define $\Sigma_i\subset {\bf P}^1 \times {\bf G}_mV$ to be the family whose
fibre over $t\in {\bf P}^1-\{0\}$ is the variety
$Y_{\alpha,i-1}(\Fdot,L_{i,t})$.

Let 
$\psi:{\bf P}^1\rightarrow \overline{U_0}={\bf G}_{n+1-m-s}V$ be a map with 
$\psi(0) = M_1: = L_{2,\infty}$,
$\psi(\infty)=K$,  and $\psi^{-1}(U_0)= {\bf P}^1-\{0\}$.
Let $\Sigma_1 \subset {\bf P}^1\times {\bf G}_mV$ be the family whose fibre
over $t\in {\bf P}^1$ is $\Omega_\alpha\Fdot\bigcap\Omega_{\psi(t)}$,
a generically transverse intersection which is irreducible for $t\neq 0$, by 
Theorem~\ref{thm:geometric_intersection}.
Then for $1\leq i\leq b$, $\Sigma_i\subset {\bf P}^1\times {\bf G}_mV$ 
is a family with equidimensional generically reduced fibres over ${\bf P}^1$.

For $1\leq i\leq b$, let 
$\phi_i:{\bf P}^1 \rightarrow \overline{{\cal G}_{i-1}}$ be the 
map associated to the family $\Sigma_i$, as in \S\ref{sec:Chow}.
Then $\phi_i(0) = \phi_{i+1}(\infty)\in {\cal G}_i$
and $\phi_i(t)\in {\cal G}_{i-1}$ for $t\neq 0$, by 
Theorem~\ref{thm:inductive_engine}.
Thus these parameterized rational curves give a chain of 
rational equivalences between 
$\Omega_\alpha\Fdot\bigcap \Omega_K$ 
and  $Y_{\alpha,b}\Fdot$.
\QED

Let $\beta\in \alpha*r$ and $\gamma\in \alpha*(r+1)$.
If $\gamma\in \beta*1$ with $j(\alpha,\gamma) = j(\beta,\gamma)$,
%equivalently, with $j(\alpha,\gamma)\leq j(\alpha,\beta)$, 
write $\beta \prec_\alpha\gamma$.

\subsection{Proof of Theorem 3.5}\label{sec:proof_of_deformations}
Let $t\neq 0$.   
Recall that $L_t$ contains the subspace 
$F_{\alpha_1 +s}$ of $\Mdot$, but not $F_{\alpha_1 +s-1}$.
Since  $M\in U_{\alpha,s-1}\Fdot$, we have 
$F_{\alpha_1}\cap M = F_{\alpha_1+s-1}$,
but $F_{\alpha_1}\cap L_t = F_{\alpha_1+s}$, thus
$F_i\cap L_t$ is a hyperplane of $F_i\cap M$ for any $i\leq \alpha_1$.
Then $L_t\in U_{\alpha,s}\Fdot$, for $t\neq 0$,  as
\begin{enumerate}
\item $F_{\alpha_1}\cap L_t=F_{\alpha_1+s}$.

\item For $1\leq j\leq m$, $F_{\alpha_j}\cap M = F_{\alpha_j+1}\cap M$.
So $F_{\alpha_j}\cap L_t = F_{\alpha_j+1}\cap L_t$.
Moreover, $\dim F_{\alpha_j}\cap L_t = \dim F_{\alpha_j}\cap M-1$,
which is $n+2-\alpha_j-j-s$.
\end{enumerate}

Suppose $t\neq 0$ and  recall that
$$
Y_{\alpha,r}(\Fdot,L_t)\ = 
\sum_{\stackrel{\mbox{ \scriptsize ${\beta}{\in}{\alpha}{*}{r}$}}
{j(\alpha,\beta)=1}}
\Omega_{\beta+ (s-1)\delta^1}\Fdot
\ +
\sum_{\stackrel{\mbox{ \scriptsize ${\beta}{\in}{\alpha}{*}{r}$}}
{j(\alpha,\beta)>1}}
X_\beta(j(\alpha,\beta),\Fdot,L_t).
$$
This defines a family $\Sigma\subset({\bf P}^1-\{0\})\times{\bf G}_mV$
with equidimensional (actually isomorphic) fibres over
${\bf P}^1-\{0\}$.
We establish Theorem~\ref{thm:inductive_engine}, showing the fibre of
$\overline{\Sigma}$ at 0 is $Y_{\alpha,r+1}(\Fdot,M)$  by examining each
component of $Y_{\alpha,r}(\Fdot,L_t)$ separately, then assembling the
result. 

Let $\beta\in \alpha*r$. 
Consider a component of   $Y_{\alpha,r}(\Fdot,L_t)$ in the first summand, so
$j(\alpha,\beta)=1$.
Then $\gamma := \beta +\delta^1$ is the unique sequence satisfying 
$\beta \prec_\alpha \gamma$.
In this case, 
$\Omega_{\beta+(s-1)\delta^1}\Fdot = 
\Omega_{\gamma+(s-2)\delta^1}\Fdot$.

Now consider a component in the second sum, so $j=j(\alpha,\beta)>1$.
Let $\beta' := \beta|_j$,
$\Fpdot := \Fdot|_{\beta_j}$, and $L'_t :=F_{\beta_j}\cap  L_t$.
For $t\neq 0$, Corollary~\ref{cor:intermediate_cycles} gives
$$
X_\beta(j(\alpha,\beta),\Fdot,L_t) \ = \ p(\pi^{-1}(
\Omega_{\beta'}\Fpdot\bigcap \Omega_{L'_t})).
$$
As $\alpha_j<\beta_j<\alpha_{j-1}$, $\dim L_t'=\dim F_{\beta_j}+1-j-(s-1)$,
by formula (\ref{sec:Y_alpha,r}).
For $1\leq i<j$, $\beta_i = \alpha_i$ and so 
$\dim L'_t\cap F_{\beta_i} =  n+2-\beta_i -i-s$.
Thus, by Theorem~\ref{thm:geometric_intersection}(3), 
$\Omega_{\beta'}\Fpdot\bigcap \Omega_{L'_t}$ is 
generically transverse and irreducible.
We study the `limit' of these cycles as $t\rightarrow 0$, in the sense of
\S\ref{sec:Chow}.
Define $L' := \lim_{t\rightarrow 0}L'_t 
= \lim_{t\rightarrow 0}  F_{\beta_j}\cap L_t$, which is 
$F_{\beta_j +1}\cap M$, by Lemma~\ref{lemma:limits_are_good}.
Then 
\begin{enumerate}
\item[(1)]  $F_{\alpha_1}\cap L' = F_{\alpha_1}\cap M = F_{\alpha_1+s-1}$.
\item[(2)]  For $1\leq i\leq j$, $F_{\beta_i}\cap L' = F_{\beta_i+1}\cap L'$.
This follows for $i=j$ because we have 
$L'\subset F_{\beta_j+1}\subset F_{\beta_j}$ and 
for $i<j$, because  $\beta_i = \alpha_i$ and 
$F_{\alpha_i}\cap M = F_{\alpha_i+1}\cap M$.
Moreover, for  $1\leq i\leq j$, 
$\dim F_{\beta_i}\cap L' =n+2-\beta_i-i-(s-1)$.
\end{enumerate}

Thus $L'\in U_{\beta',s-1}\Fpdot$ so
$\Omega_{\beta'}\Fpdot \bigcap \Omega_{L'}$ is generically 
transverse, by Theorem~\ref{thm:geometric_intersection}(1).
So, 
$$
\lim_{t\rightarrow 0} X_\beta(j(\alpha,\beta),\Fdot,L_t) \ =\  p(\pi^{-1}(
\Omega_{\beta'}\Fpdot\bigcap \Omega_{L'})).
$$
But $\Span{F_{\beta_{i-1}},F_{\beta_i}\cap L}
\subset F_{\beta_i +1}$, since $L'\in U_{\beta',s-1}\Fpdot$.
By Remark~\ref{remarkI},
$$
\Omega_{\beta'}\Fpdot\bigcap \Omega_{L'}
\ =\  \sum_{\stackrel{\mbox{\scriptsize 
${\gamma'}{\in}{\beta'}{*}{1}$}}{j(\beta',\gamma')=1}}
\Omega_{\gamma'+(s-2)\delta^1}\Fdot \ + 
\sum_{\stackrel{\mbox{\scriptsize 
${\gamma'}{\in}{\beta'}{*}{1}$}}{j(\beta',\gamma')>1}}
X_{\gamma'}(j(\beta',\gamma'),\Fpdot,L').
$$
And so 
$
\lim_{t\rightarrow 0} X_\beta(j(\alpha,\beta),\Fdot,L_t)$
is the cycle
$$
\sum_{\stackrel{\mbox{\scriptsize 
${\gamma'}{\in}{\beta'}{*}{1}$}}{j(\beta',\gamma')=1}}
p(\pi^{-1}(\Omega_{\gamma'+(s-2)\delta^1}\Fpdot)) \ + 
\sum_{\stackrel{\mbox{\scriptsize 
${\gamma'}{\in}{\beta'}{*}{1}$}}{j(\beta',\gamma')>1}}
p(\pi^{-1}(X_{\gamma'}(j(\beta',\gamma'),\Fpdot,L'))).
$$

We simplify this expression, beginning with the first sum.
Let $\gamma'\in \beta'*1$ satisfy $j(\beta',\gamma')=1$.
Then by Lemma~\ref{lemma:rational_fibration}, 
$p(\pi^{-1}(\Omega_{\gamma'+(s-2)\delta^1}\Fdot))$ equals
$ \Omega_{\gamma+(s-2)\delta^1}\Fdot$,
where $\gamma:=\beta+\delta^1$ is the unique sequence 
with $\beta\prec_\alpha\gamma$ and 
$j(\alpha,\gamma)=1$.

Consider terms in the second sum,  those for which 
$\gamma'\in \beta'*1$ with $j(\beta',\gamma')>1$.   
Then $p(\pi^{-1}(X_{\gamma'}(j(\beta',\gamma'),\Fpdot,L')))$ is 
the subvariety of $\Omega_\beta\Fdot$ consisting of those 
$H$ such that 
there exists $ K\subset H$ with $\dim K=j$,
$K\in \Omega_{\gamma'}\Fpdot$, and 
$\dim K\cap F'_{\gamma'_{j(\beta',\gamma')}}\cap L' \geq 1$.
\smallskip

Let $\gamma := \beta+\delta^{j(\beta',\gamma')}$,
the unique sequence with $\beta\prec_\alpha \gamma$ 
and $j(\alpha,\gamma)=j(\beta',\gamma')$.
Then, as $\gamma_{j(\alpha,\gamma)}>\beta_j$, 
the definition of $\Fpdot$ implies 
$F'_{\gamma'_{j(\beta',\gamma')}}=F_{j(\alpha,\gamma)}
\subset F_{\beta_j+1}$.
Since $L'=F_{\beta_j+1}\cap M$, we see that 
$F'_{\gamma'_{j(\beta',\gamma')}}\cap L' = 
F_{\gamma_{j(\alpha,\gamma)}}\cap M$.
Thus if 
$$
H\ \in\ p(\pi^{-1}(X_{\gamma'}(j(\beta',\gamma'),\Fpdot,L'))),
$$
then $H\in \Omega_\gamma \Fdot$ and 
$\dim H\cap F_{\gamma_{j(\alpha,\gamma)}}\cap M\geq 1$, so 
$H\in  X_\gamma(j(\alpha,\gamma), \Fdot,M)$.
The reverse inclusion, 
$$
X_\gamma(j(\alpha,\gamma), \Fdot,M)\ \subset\  
p(\pi^{-1}(X_{\gamma'}(j(\beta',\gamma'),\Fpdot,L'))),
$$
is similar.

This shows that 
$\lim_{t\rightarrow 0} X_{\beta}(j(\alpha,\beta),\Fdot,L_t)$
is the cycle 
$$
(\ref{sec:proof_of_deformations})\hspace{1.1in}
\sum_{\stackrel{\mbox{\scriptsize${\beta}{\prec_\alpha}{\gamma}$}}
{j(\alpha,\gamma)=1}} 
\Omega_{\gamma+ (s-2)\delta^1}\Fdot
\ +
\sum_{\stackrel{\mbox{\scriptsize${\beta}{\prec_\alpha}{\gamma}$}}
{j(\alpha,\gamma)>1}} 
X_\gamma(j(\alpha,\gamma),\Fdot,L).
\hspace{1.4in}
$$
The sets $\{\gamma \,|\, \beta\prec_\alpha\gamma\}$ for 
$\beta\in \alpha*r$ partition
the set $\alpha*(r+1)$.
Thus 
$$
\lim_{t\rightarrow 0}
Y_{\alpha,r}(\Fdot,L_t)\ = 
\sum_{\stackrel{\mbox{ \scriptsize ${\gamma}{\in}{\alpha}{*}{(r{+}1)}$}}
{j(\alpha,\gamma)=1}}
\Omega_{\gamma+ (s-2)\delta^1}\Fdot
\ +
\sum_{\stackrel{\mbox{ \scriptsize ${\gamma}{\in}{\alpha}{*}{(r{+}1)}$}}
{j(\alpha,\gamma)>1}}
X_\beta(j(\alpha,\gamma),\Fdot,M),
$$
which is $Y_{\alpha,r+1}(\Fdot,M)$.
\QED

\section{Link to Schensted insertion}\label{sec:schensted}
The set ${[n]\choose m}$ has a partial order, called the Bruhat order:
$\alpha\leq \beta$ if and only if 
$\Omega_\beta\Fdot \subset \Omega_\alpha\Fdot$.
Combinatorially, this is $\alpha\leq \beta$ if $\alpha_i\leq \beta_i$ for 
$1\leq i\leq m$.

We interpret the behavior of the components 
$X_\beta(j(\alpha,\beta)\Fdot,L)$ of 
the intermediate varieties $Y_{\alpha,i-1}(\Fdot,L)$ 
in our proof of Pieri's formula (\S\ref{sec:pieri_proof}) 
as the branching of a certain subtree of 
${[n]\choose m}$ with root $\alpha$.
This tree arises similarly in the combinatorial proof of 
Pieri's formula for Schur polynomials using Schensted 
insertion given in~\cite{Fulton_tableaux}.
To simplify this discussion, assume further that $n>\alpha_1+b$.

Each rational equivalence of
\S\ref{sec:pieri_proof} is induced by a family $\Sigma_i$ over ${\bf P}^1$
with generic fibre  in ${\cal G}_{i-1}$ and special 
fibre in ${\cal G}_i$.
The components of cycles in ${\cal G}_{i-1}$ are indexed by 
$\beta\in \alpha*(i-1)$, with 
$\beta$th component  $\Omega_{\beta+(b+1-i)\delta^1}\Fdot$, if 
$j(\alpha,\beta)=1$, and $X_{\beta}(j(\alpha,\beta),\Fdot,L)$ 
otherwise.
In passing to ${\cal G}_i$ via $\phi_{i}$, the component 
$\Omega_{\beta+(b+1-i)\delta^1}\Fdot$ is unchanged, but reindexed:
$\Omega_{\gamma+(b-i)\delta^1}\Fdot$, where 
$\gamma := \beta+\delta^1$ is the unique sequence with 
$\beta\prec_\alpha\gamma$.
By equation (\ref{sec:proof_of_deformations}), the other components
become 
$$
\sum_{\stackrel{\mbox{\scriptsize$\beta{\prec_\alpha}\gamma$}}
{j(\alpha,\gamma)=1}}
\Omega_{\gamma+ (b-i)\delta^1}\Fdot \ + 
\sum_{\stackrel{\mbox{\scriptsize$\beta{\prec_\alpha}\gamma$}}
{j(\alpha,\gamma)>1}}X_\gamma(j(\alpha,\gamma),\Fdot,M_i).
$$
Thus the component of the generic fibre of $\Sigma_i$
indexed by $\beta\in \alpha*(i-1)$  
becomes a sum of 
components indexed by $\{\gamma\,|\, \beta\prec_\alpha\gamma\}$ 
at the special fibre.
\smallskip

This suggest defining a tree  ${\cal T}_{\alpha,b}$ whose branching
represents the `branching' of components of $Y_{\alpha,i-1}(\Fdot,L)$  in
these deformations.
Let ${\cal T}_{\alpha,b}\subset {[n]\choose m}$ be the tree with
vertex set $\bigcup\{\alpha*i\,|\,0\leq i\leq b\}$ and covering relation
$\beta\prec_\alpha\gamma$.
This is a tree as $\alpha*i$ is partitioned by the sets 
$\{\gamma \,|\, \beta\prec_\alpha\gamma\}$ for $\beta\in \alpha*(i-1)$.

For a decreasing $m$-sequence $\alpha$, let $\lambda(\alpha)$ be the partition
$(\alpha_1-m,\alpha_2-m+1,\ldots,\alpha_m-1)$.
The association $\alpha\longleftrightarrow \lambda(\alpha)$ 
gives an order isomorphism between the set of decreasing 
$m$-sequences and the set of  partitions of length at most $m$.
This transfers  notions for sequences into 
corresponding notions for partitions.

To a (semi-standard) Young tableau $T$ with entries among $1,\ldots,m$, 
associate a monomial $x^T$ in the variables
$x_1,x_2\ldots,x_m$:
The exponent of $x_i$ in $x^T$ is the number of occurrences of 
$i$ in $T$.
The Schur polynomial $s_\lambda$ is 
$\sum x^T$, the sum over all tableaux $T$ of shape $\lambda$.
There is surjective homomorphism from the algebra of Schur polynomials
to the Chow ring of ${\bf G}_mV$ defined by:
$$
s_\lambda \longmapsto \left\{\begin{array}{ll} 
[\Omega_\alpha\Fdot]&\mbox{ \ if \ } \lambda = \lambda(\alpha)
\mbox{ for some }\alpha\in {[n]\choose m}\\
0&\mbox{ \ otherwise}\end{array}\right..
$$
Special Schur polynomials are indexed by partitions $(b,0,\ldots,0)$ with 
a single row.

Schensted insertion gives a combinatorial proof of Pieri's 
formula, providing a content-preserving bijection between the 
set of pairs $(S,T)$ of tableaux where $S$ has shape $\lambda$ and 
$T$ has shape $(b,0\ldots,0)$ and the set of all tableaux whose shape is in 
$\lambda*b\,$:
Insert the reading word of $T$
into $S$.
The resulting tableau has shape $\mu\in \lambda*b$.

Let $\lambda = \lambda_0,\lambda_1,\ldots,\lambda_b=\mu$ be the sequence
of shapes resulting from the  
insertion of successive entries of $T$ into $S$.
Since $T$ is a single row, it is a property
of the insertion algorithm that 
$\lambda_i\prec_\lambda \lambda_{i+1}$,
and so this sequence is a chain in the tree ${\cal T}_{\alpha,b}$.

The totality of these insertions for all such pairs of tableaux 
gives all chains in  ${\cal T}_{\alpha,b}$.
Thus the `branching' of shapes during Schensted insertion 
is identical to the branching of components 
in the rational equivalences of \S\ref{sec:pieri_proof}.
We feel this relation to combinatorics is one of the more intriguing 
aspects of our proof of Pieri's formula.
It leads us to speculate that similar ideas may yield a 
geometric proof of the Littlewood-Richardson rule.

%\bibliographystyle{siam}
%\bibliography{/home/sottile/papers/bibl/bibliography,%
%/home/sottile/papers/bibl/sottile}
\vfill
\pagebreak

\begin{center}
\sc Appendices
\end{center}

Following are two appendices which are for general distribution and not
intended for publication.
In the first, we illustrate how one may introduce coordinates by giving an
explicit parameterized family of intersections in one example and studying
how the intersection cycles deform.
It should be possible to transfer this description into a parameterized
family of ideals, and use computer algebra to study the deformations.

In the second, we illustrate how deformations of the kind studied here may be
used to find real solutions to some problems of enumerative geometry.

\begin{center}
\sc Appendix A: An example illustrating explicit rational equivalence
\end{center}

Our proof of Pieri's formula used methods which 
were  almost\footnote{
Coordinates were only chosen in the proof of 
Lemma~3.4.
A synthetic proof, while possible, is perhaps less convincing, 
and certainly not as brief.}
entirely synthetic (coordinate-free).
Here,  we illustrate our methods on a particular example,
simultaneously showing how to introduce coordinates.
While we do not attempt such an analysis, it should be possible to
replicate our results using (computer) algebra.
It is in this context that these `explicit rational equivalences' are 
most explicit.
\bigskip

\noindent{\bf A.1 Example: \boldmath $n=9$, $m=3$, $s=2$, and $\alpha = 741$.}
Let $\Fdot$ be a flag and $L$ a five-dimensional subspace of 
$V = k^9$ which meets $\Fdot$ properly.
We choose coordinates for $V$ and give a two parameter family of 
5-planes containing $L$ which exhibits a `chain of rational equivalences'
between $\Omega_{741}\Fdot\bigcap \Omega_L$ and
$Y_{741,\,2}\Fdot = \sum_{\beta\in 741*2}\Omega_\beta\Fdot$, 
which is the cycle 
$$
\Omega_{941}\Fdot\  +\ \Omega_{851}\Fdot\ +\ \Omega_{761}\Fdot\ +\
\Omega_{842}\Fdot\ +\ \Omega_{752}\Fdot\ +\ \Omega_{743}\Fdot.
$$

Let $e_1,\ldots,e_9$ be a basis for $V$ such that 
$F_j = \Span{e_j,\ldots,e_9}$, for $1\leq j\leq 9$.
Let $x_1,\ldots,x_9$ be dual coordinates for $V^*$, then 
$F_j$ is defined by the  linear forms $x_1,\ldots, x_{j-1}$.
Further assume that $e_1,\ldots, e_9$ have been chosen so that
$L = \Span{e_1,e_2,e_3,e_4,e_5}$.

For $(s,t)\in k^2$,  define 
linearly independent forms $\Lambda_1,\Lambda_2,\Lambda_3$, and $\Lambda_4$
in $V^*$:
\begin{eqnarray*}
\Lambda_1 &:=& x_1 + s x_8\\
\Lambda_2 &:=& x_2 + t x_3 + st^2x_4+(t^2+st^3)x_5+(t^3+st^4)x_6+t^4 x_8\\
\Lambda_3 &:=& x_4 + s x_9\\
\Lambda_4 &:=& x_7
\end{eqnarray*}
If $s\cdot t\neq 0$, let $\sigma=1/s$ and $\tau = 1/t$.
Consider the forms
\begin{eqnarray*}
\sigma\tau^4\Lambda_2 - \sigma^2\Lambda_1 &=& -\sigma^2 x_1 +\sigma\tau^4 x_2 
+ \sigma\tau^3 x_3 + \tau^2x_4 +
(\sigma\tau^2+\tau)x_5 + (\sigma\tau+1)x_6\\
\Lambda_4 &=& x_7\\
\sigma\Lambda_1 &=& \sigma x_1 + x_8\\
\sigma\Lambda_3 &=& \sigma x_4 +  x_9.
\end{eqnarray*}
These are linearly independent for  $(\sigma,\tau)\in k^2$.
These two  sets of forms together  define a family $L_{s,t}$
of 5-dimensional subspaces of $k^9$ with base 
$U:= {\bf P}^1\times{\bf P}^1 - \{(\infty,0),(0,\infty)\}$, and with 
$L_{\infty,\infty} = L$.
For $s\cdot t\neq 0$,
the second set of forms shows that $L_{s,t}$ meets 
the subspaces $F_1, F_4$, and $F_7$ properly.
This shows that for $s\cdot t\neq 0$ all 5-planes $L_{s,t}$
are in the dense orbit of the parabolic subgroup $P$ stabilizing 
 $F_1, F_4$, and $F_7$.

Thus over $({\bf P}^1-\{0\})\times({\bf P}^1-\{0\})$, 
there is a family $\Sigma$ 
of cycles with fibre over $(s,t)$:
$$
\hspace{2.2in}
\Sigma_{s,t}\ :=\ \Omega_{741}\Fdot \bigcap \Omega_{L_{s,t}}.
\hspace{1.9in}(A.1)
$$
Moreover, all fibres of $\Sigma$ are isomorphic by some element of $P$.
We study the cycles obtained as first $s\rightarrow 0$ with 
$t\neq 0$ fixed, and then $t\rightarrow 0$.
We  show:
\begin{enumerate}
\item[(A)]  For $s\cdot t\neq 0$, the  intersection
$\Omega_{741}\Fdot \bigcap \Omega_{L_{s,t}}$
is generically transverse and irreducible.
\item[(B)]  For $t\neq 0$ or $\infty$, the intersection cycle 
$Y_t := \Omega_{741}\Fdot \bigcap \Omega_{L_{0,t}}$ 
is generically transverse, but 
reducible, with components indexed by $\{841, 751, 742\}=741 *1$.
This is  because if $t\neq 0$ or $\infty$, then $L_{0,t} \in U_{741,2}\Fdot$.
Thus $\Sigma$ can be extended to  $U-\{k\times 0\}$ by 
(\S A.1).
\item[(C)] The intersection $\Omega_{741}\Fdot \bigcap \Omega_{L_{0,0}}$
is improper.   
However, $\lim_{t\rightarrow 0} Y_t=
\sum_{\beta\in 741*2}\Omega_\beta\Fdot$.
\end{enumerate}
We show these points:
\medskip 

\noindent{\bf (A):}
We noted already that for $s\cdot t\neq 0$, $L_{s,t}$ meets 
$F_1, F_4$, and $F_7$ properly.
Then by Theorem~2.3, the intersection
$\Omega_{741}\Fdot \bigcap \Omega_{L_{s,t}}$ is generically transverse 
and irreducible.
\medskip 

\noindent{\bf (B):}
Fix $t\neq 0,\infty$, and consider the family $L_{s,t}$ for 
$s\in {\bf P}^1$.
Then $L_{0,t}$ is given by the forms:
$$
x_1,\ \ x_4,\ \ x_7,\ \ x_2+tx_3+t^2x_5+t^3x_6+t^4x_8
$$
and so it does not meet $\Fdot$ properly.
Note that $L_{0,t}$ has a basis:
$$
te_2-e_3,\ \ te_3-e_5,\ \  te_5-e_6,\ \ te_6-e_8,\ \ e_9.
$$
Then $L_{0,t}\in U_{741,\,2}\Fdot$ as 
\begin{enumerate}
\item $L_{0,t}\subset F_2$.
\item $L_{0,t}\cap F_4 = L_{0,t}\cap F_5 = 
\Span{te_5-e_6,\ \ te_6-e_8,\ \ e_9}$,
which has dimension 3.
\item $L_{0,t}\cap F_7 = F_9$.
\end{enumerate}
Note also that $F_2 = \Span{L_{0,t},F_4}$, 
$F_5 = \Span{L_{0,t}\cap F_4,F_7}$, and 
$F_8 \supset \Span{L_{0,t}\cap F_7}$.
By Theorem~2.3,
$$
\Omega_{741}\Fdot \bigcap \Omega_{L_{0,t}} \ =\ 
X_{841}(1,\Fdot,L_{0,t}) + X_{751}(2,\Fdot,L_{0,t}) +
X_{742}(3,\Fdot,L_{0,t}).
$$

\noindent{\bf (C):}
Note  $L_{0,0}=\Span{e_3,e_5,e_6,\,e_8,e_9}$.
Since $F_8\subset L_{0,0}$, 
$\Omega_{841}\Fdot\subset \Omega_{741}\Fdot \bigcap \Omega_{L_{0,0}}$, 
the intersection $\Omega_{741}\Fdot \bigcap \Omega_{L_{0,0}}$ is improper.

Instead let $\Sigma$ be the family of cycles over 
${\bf P}^1$ 
whose fibre over $t\neq 0,\infty$ is 
$$
\Sigma_t\ :=\ \Omega_{741}\Fdot \bigcap \Omega_{L_{0,t}},
$$
and we describe the fibre $\Sigma_0\subset 
\Omega_{741}\Fdot \bigcap \Omega_{L_{0,0}}$.
For this,  we need a better description of the fibres
$\Sigma_t$, for $t\neq 0,\infty$.
First, let $M = \Span{e_2,e_3,e_5,e_6, e_8, e_9}\in U_{741,1}\Fdot$
\vspace{12pt}

\begin{enumerate}
\item[$j=1$:] 
$X_{841}(1,\Fdot,L_{0,t}) = \{H\in \Omega_{841}\Fdot\,|\, 
\dim H\cap F_8\cap L_{0,t}\geq 1\}$.
Since $L_{0,t}\cap F_8 = F_9$, this equals 
$\Omega_{941}\Fdot = \Omega_{841+100}\Fdot$.
\item[$j=2$:] Similarly, 
$X_{751}(2,\Fdot,L_{0,t}) = \{H\in \Omega_{751}\Fdot\,|\, 
\dim H\cap F_5\cap L_{0,t}\geq 1\}$.
If $p,\pi$ are the first and second projections on
$$
\widetilde{\Omega_{751}}^2\Fdot\ =\ 
\{ (H,K)\in\Omega_{751}\Fdot\times{\bf G}_3F_5
\,|\, K\subset H, \  \dim K\cap F_7 \geq 1\},
$$
then $X_{751}(2,\Fdot,L_{0,t}) = p(\pi^{-1}(\Omega_{31}\Fdot|_5
\bigcap \Omega_{F_5\cap L_{0,t}}))$.

Since $F_5\cap L_{0,t} = \Span{e_6-te_5,e_8-te_6, e_9}$,
$\lim_{t\rightarrow 0}F_5\cap L_{0,t} = F_6\cap M$.

Suppose $K\in \Omega_{31}\Fdot|_5\bigcap \Omega_{F_6\cap M}$.
Then $K\subset F_5$, $\dim K\cap F_7\geq 1$, and 
$\dim K\cap F_6\cap M\geq 1$.
\begin{enumerate}
\item If $K\cap (F_6\cap M) \subset F_7$, then $\dim K\cap F_8\geq 1$,
as $F_6\cap M\cap F_7 = F_8$, 
and so $K\in \Omega_{41}\Fdot|_5$.
Also,  $p(\pi^{-1}(\Omega_{41}\Fdot|_5)) = \Omega_{851}\Fdot$.
\item If $K\cap F_6\cap M \not\subset F_7$, then 
$K\subset \Span{F_7, M\cap F_6} = F_6$ and so 
$K\in \Omega_{32}\Fdot|_5$.
Note that $p(\pi^{-1}(\Omega_{32}\Fdot|_5)) = \Omega_{761}\Fdot$.
\end{enumerate}
This shows that $\Omega_{31}\Fdot|_5\bigcap \Omega_{F_6\cap M}
= \Omega_{41}\Fdot|_5+\Omega_{32}\Fdot|_5$.
Since this intersection is generically transverse, it follows that 
$$
\lim_{t\rightarrow 0}  X_{751}(2,\Fdot,L_{0,t})= 
p(\pi^{-1}(\Omega_{41}\Fdot|_5+\Omega_{21}\Fdot|_5))
\ =\  \Omega_{851}\Fdot+\Omega_{761}\Fdot.
$$
\item[$j=3$:] 
$X_{742}(3,\Fdot,L_{0,t}) = \{H\in\Omega_{742}\Fdot\,|\, 
\dim H\cap L_{0,t}\cap F_2 \geq 1\}$.
This is $\Omega_{631}\Fdot|_2\bigcap \Omega_{L_{0,t}}$, which is 
generically transverse and irreducible, {\em as an 
intersection in ${\bf G}_3F_2$}, by 
Theorem~2.3.

Similarly,  $X_{742}(3,\Fdot,L_{0,0})=
\Omega_{631}\Fdot|_2\bigcap \Omega_{L_{0,0}}$ is 
generically transverse but reducible, as an 
intersection in ${\bf G}_3F_2$,
by Theorem~2.3(2).
Thus 
$$
\lim_{t\rightarrow 0}X_{742}(3,\Fdot,L_{0,t}) \ =\ 
X_{742}(3,\Fdot,L_{0,0}).
$$
Let $H\in X_{742}(3,\Fdot,L_{0,0})$.
Then $H\subset F_2$, $\dim H\cap F_4\geq 2$, $\dim H\cap F_7\geq 1$, and 
$\dim H\cap L_{0,0}\geq 1$.
\begin{enumerate}
\item If $H\cap L_{0,0}\subset F_7$, then $\dim H\cap F_8\geq 1$, as 
$F_7\cap L_{0,0}=F_8$.
Thus $H\in \Omega_{842}\Fdot$.
\item If $\dim H\cap F_4\cap  L_{0,0} \geq 1$ and 
$H\cap  L_{0,0}\not\subset F_7$, then 
$\dim H\cap \Span{F_7,F_4\cap L_{0,0}}\geq 2$.    
Since $F_5 = \Span{F_7,F_4\cap L_{0,0}}$, 
$H\in \Omega_{752}\Fdot$.
\item If $H\cap L_{0,0}\not\subset F_4$, then $H\subset 
\Span{F_4,F_2\cap L_{0,0}} = F_3$,
and so $H\in \Omega_{743}\Fdot$.
\end{enumerate}
We deduce that $X_{742}(3,\Fdot,L_{0,0})
= \Omega_{842}\Fdot+ \Omega_{752}\Fdot + \Omega_{743}\Fdot$, and so 
$$
\lim_{t\rightarrow 0}  X_{742}(3,\Fdot,L_{0,t})\ = \
 \Omega_{842}\Fdot+ \Omega_{752}\Fdot + \Omega_{743}\Fdot.
$$
\end{enumerate}

These calculations identify the fibre $\Sigma_0$ of the family
whose fibre over $t\neq 0$ is 
$\Omega_\alpha\Fdot\bigcap \Omega_{L_{0,t}}$.
They show $\Sigma_0 =  \Omega_{841+100}\Fdot + 
X_{751}(2,\Fdot,L_{0,0})+X_{742}(3,\Fdot,L_{0,0})$, which is 
$$
\Omega_{941}\Fdot\ +\ \Omega_{851}\Fdot\ +\ \Omega_{761}\Fdot\ +\ 
\Omega_{842}\Fdot\ +\ \Omega_{752}\Fdot\ +\ \Omega_{743}\Fdot\ \
=\ \sum_{\beta\in 741*2}\Omega_\beta\Fdot.
$$

\noindent{\bf A.2: The rational equivalences, explicitly}
The calculations of \S A.1
yield an explicit chain of 
rational equivalences from $\Omega_{741}\Fdot\bigcap \Omega_L$
to $\sum_{\beta\in 741*2}\Omega_\beta\Fdot$.
Furthermore, the generic fibre of $\Sigma_1$ is irreducible 
and the special fibre has three components, indexed by the set $741*1$.
The generic fibre of $\Sigma_2$ also has three components indexed by 
$741*1$, and its special fibre has six
components, indexed by the set $741*2$.
The tree ${\cal T}_{741,2}$, shown in Figure 1 in the lattice of 
Young diagrams, represents this branching of components.
\begin{figure}[htb]
$$
\setlength{\unitlength}{0.008in}%
\begin{picture}(760,404)(100,396)
\thicklines
 \put(400,730){\line(-1,-1){85}}
 \put(475,650){\line(-1, 2){38}}
 \put(610,650){\line(-5, 3){150}}
 \put(150,470){\line( 1, 1){100}}
 \put(270,570){\line(0,-1){100}}
 \put(288,572){\line( 3,-5){62}}
 \put(480,590){\line( 1,-5){24}}
 \put(620,470){\line(-5, 6){105}}
 \put(760,470){\line(-5, 6){105}}
\put(260,600){\line( 1, 0){ 20}}
\put(280,600){\line( 0,-1){ 20}}
\put(280,580){\line(-1, 0){ 20}}
\put(260,580){\line( 0, 1){ 20}}
\put(100,420){\line( 1, 0){ 20}}
\put(120,420){\line( 0,-1){ 20}}
\put(120,400){\line(-1, 0){ 20}}
\put(100,400){\line( 0, 1){ 20}}
\put(480,620){\line( 1, 0){ 20}}
\put(500,620){\line( 0,-1){ 20}}
\put(500,600){\line(-1, 0){ 20}}
\put(480,600){\line( 0, 1){ 20}}
\put(680,640){\line( 1, 0){ 20}}
\put(700,640){\line( 0,-1){ 20}}
\put(700,620){\line(-1, 0){ 20}}
\put(680,620){\line( 0, 1){ 20}}
\put(840,460){\line( 1, 0){ 20}}
\put(860,460){\line( 0,-1){ 20}}
\put(860,440){\line(-1, 0){ 20}}
\put(840,440){\line( 0, 1){ 20}}
\put(820,460){\line( 1, 0){ 20}}
\put(840,460){\line( 0,-1){ 20}}
\put(840,440){\line(-1, 0){ 20}}
\put(820,440){\line( 0, 1){ 20}}
\put(680,460){\line( 1, 0){ 20}}
\put(700,460){\line( 0,-1){ 20}}
\put(700,440){\line(-1, 0){ 20}}
\put(680,440){\line( 0, 1){ 20}}
\put(640,440){\line( 1, 0){ 20}}
\put(660,440){\line( 0,-1){ 20}}
\put(660,420){\line(-1, 0){ 20}}
\put(640,420){\line( 0, 1){ 20}}
\put(540,440){\line( 1, 0){ 20}}
\put(560,440){\line( 0,-1){ 20}}
\put(560,420){\line(-1, 0){ 20}}
\put(540,420){\line( 0, 1){ 20}}
\put(520,440){\line( 1, 0){ 20}}
\put(540,440){\line( 0,-1){ 20}}
\put(540,420){\line(-1, 0){ 20}}
\put(520,420){\line( 0, 1){ 20}}
\put(420,460){\line( 1, 0){ 20}}
\put(440,460){\line( 0,-1){ 20}}

\put(440,440){\line(-1, 0){ 20}}
\put(420,440){\line( 0, 1){ 20}}
\put(340,420){\line( 1, 0){ 20}}
\put(360,420){\line( 0,-1){ 20}}
\put(360,400){\line(-1, 0){ 20}}
\put(340,400){\line( 0, 1){ 20}}
\put(260,440){\line( 1, 0){ 20}}
\put(280,440){\line( 0,-1){ 20}}
\put(280,420){\line(-1, 0){ 20}}
\put(260,420){\line( 0, 1){ 20}}
\put(220,420){\line( 1, 0){ 20}}
\put(240,420){\line( 0,-1){ 20}}
\put(240,400){\line(-1, 0){ 20}}
\put(220,400){\line( 0, 1){ 20}}
\put(120,420){\line( 1, 0){ 20}}
\put(140,420){\line( 0,-1){ 20}}
\put(140,400){\line(-1, 0){ 20}}
\put(120,400){\line( 0, 1){ 20}}
\put(620,460){\line( 0,-1){ 40}}
\put(740,440){\line( 1, 0){ 80}}
\put(820,440){\line( 0, 1){ 20}}
\put(820,460){\line(-1, 0){ 80}}
\put(740,460){\line( 0,-1){ 40}}
\put(740,420){\line( 1, 0){ 40}}
\put(780,420){\line( 0, 1){ 40}}
\put(800,460){\line( 0,-1){ 20}}
\put(760,460){\line( 0,-1){ 40}}
\put(600,620){\line( 1, 0){ 80}}
\put(680,620){\line( 0, 1){ 20}}
\put(680,640){\line(-1, 0){ 80}}
\put(600,640){\line( 0,-1){ 40}}
\put(600,600){\line( 1, 0){ 40}}
\put(640,600){\line( 0, 1){ 40}}
\put(660,640){\line( 0,-1){ 20}}
\put(400,760){\line( 1, 0){ 80}}
\put(480,760){\line( 0, 1){ 20}}
\put(480,780){\line(-1, 0){ 80}}
\put(400,780){\line( 0,-1){ 40}}
\put(400,740){\line( 1, 0){ 40}}
\put(440,740){\line( 0, 1){ 40}}
\put(620,640){\line( 0,-1){ 40}}
\put(500,640){\line( 0,-1){ 20}}
\put(460,640){\line( 0,-1){ 40}}
\put(260,620){\line( 1, 0){ 80}}
\put(340,620){\line( 0, 1){ 20}}
\put(340,640){\line(-1, 0){ 80}}
\put(260,640){\line( 0,-1){ 40}}
\put(260,600){\line( 1, 0){ 40}}
\put(300,600){\line( 0, 1){ 40}}
\put(320,640){\line( 0,-1){ 20}}
\put(280,640){\line( 0,-1){ 40}}
\put(100,440){\line( 1, 0){ 80}}
\put(180,440){\line( 0, 1){ 20}}
\put(180,460){\line(-1, 0){ 80}}
\put(100,460){\line( 0,-1){ 40}}
\put(100,420){\line( 1, 0){ 40}}
\put(140,420){\line( 0, 1){ 40}}
\put(440,620){\line( 1, 0){ 80}}
\put(520,620){\line( 0, 1){ 20}}
\put(520,640){\line(-1, 0){ 80}}
\put(440,640){\line( 0,-1){ 40}}

\put(440,600){\line( 1, 0){ 40}}
\put(480,600){\line( 0, 1){ 40}}
\put(420,780){\line( 0,-1){ 40}}
\put(460,780){\line( 0,-1){ 20}}
\put(160,460){\line( 0,-1){ 20}}
\put(120,460){\line( 0,-1){ 40}}
\put(220,440){\line( 1, 0){ 80}}
\put(300,440){\line( 0, 1){ 20}}
\put(300,460){\line(-1, 0){ 80}}
\put(220,460){\line( 0,-1){ 40}}
\put(220,420){\line( 1, 0){ 40}}
\put(260,420){\line( 0, 1){ 40}}
\put(280,460){\line( 0,-1){ 20}}
\put(240,460){\line( 0,-1){ 40}}
\put(340,440){\line( 1, 0){ 80}}
\put(420,440){\line( 0, 1){ 20}}
\put(420,460){\line(-1, 0){ 80}}
\put(340,460){\line( 0,-1){ 40}}
\put(340,420){\line( 1, 0){ 40}}
\put(380,420){\line( 0, 1){ 40}}
\put(400,460){\line( 0,-1){ 20}}
\put(360,460){\line( 0,-1){ 40}}
\put(480,440){\line( 1, 0){ 80}}
\put(560,440){\line( 0, 1){ 20}}
\put(560,460){\line(-1, 0){ 80}}
\put(480,460){\line( 0,-1){ 40}}
\put(480,420){\line( 1, 0){ 40}}
\put(520,420){\line( 0, 1){ 40}}
\put(540,460){\line( 0,-1){ 20}}
\put(500,460){\line( 0,-1){ 40}}
\put(600,440){\line( 1, 0){ 80}}
\put(680,440){\line( 0, 1){ 20}}
\put(680,460){\line(-1, 0){ 80}}
\put(600,460){\line( 0,-1){ 40}}
\put(600,420){\line( 1, 0){ 40}}
\put(640,420){\line( 0, 1){ 40}}
\put(660,460){\line( 0,-1){ 20}}
\put(846,444){\makebox(0,0)[lb]{\raisebox{0pt}[0pt][0pt]{\small\bf 2}}}
\put(105,404){\makebox(0,0)[lb]{\raisebox{0pt}[0pt][0pt]{\small\bf 1}}}
\put(345,404){\makebox(0,0)[lb]{\raisebox{0pt}[0pt][0pt]{\small\bf 1}}}
\put(525,424){\makebox(0,0)[lb]{\raisebox{0pt}[0pt][0pt]{\small\bf 1}}}
\put(645,424){\makebox(0,0)[lb]{\raisebox{0pt}[0pt][0pt]{\small\bf 1}}}
\put(825,444){\makebox(0,0)[lb]{\raisebox{0pt}[0pt][0pt]{\small\bf 1}}}
\put(265,584){\makebox(0,0)[lb]{\raisebox{0pt}[0pt][0pt]{\small\bf 1}}}
\put(485,604){\makebox(0,0)[lb]{\raisebox{0pt}[0pt][0pt]{\small\bf 1}}}
\put(685,624){\makebox(0,0)[lb]{\raisebox{0pt}[0pt][0pt]{\small\bf 1}}}
\put(225,404){\makebox(0,0)[lb]{\raisebox{0pt}[0pt][0pt]{\small\bf 1}}}
\put(125,404){\makebox(0,0)[lb]{\raisebox{0pt}[0pt][0pt]{\small\bf 2}}}
\put(265,424){\makebox(0,0)[lb]{\raisebox{0pt}[0pt][0pt]{\small\bf 2}}}
\put(425,444){\makebox(0,0)[lb]{\raisebox{0pt}[0pt][0pt]{\small\bf 2}}}
\put(545,424){\makebox(0,0)[lb]{\raisebox{0pt}[0pt][0pt]{\small\bf 2}}}
\put(685,444){\makebox(0,0)[lb]{\raisebox{0pt}[0pt][0pt]{\small\bf 2}}}
\end{picture}
$$
\caption{The Tree ${\cal T}_{741,2}$}
\end{figure}
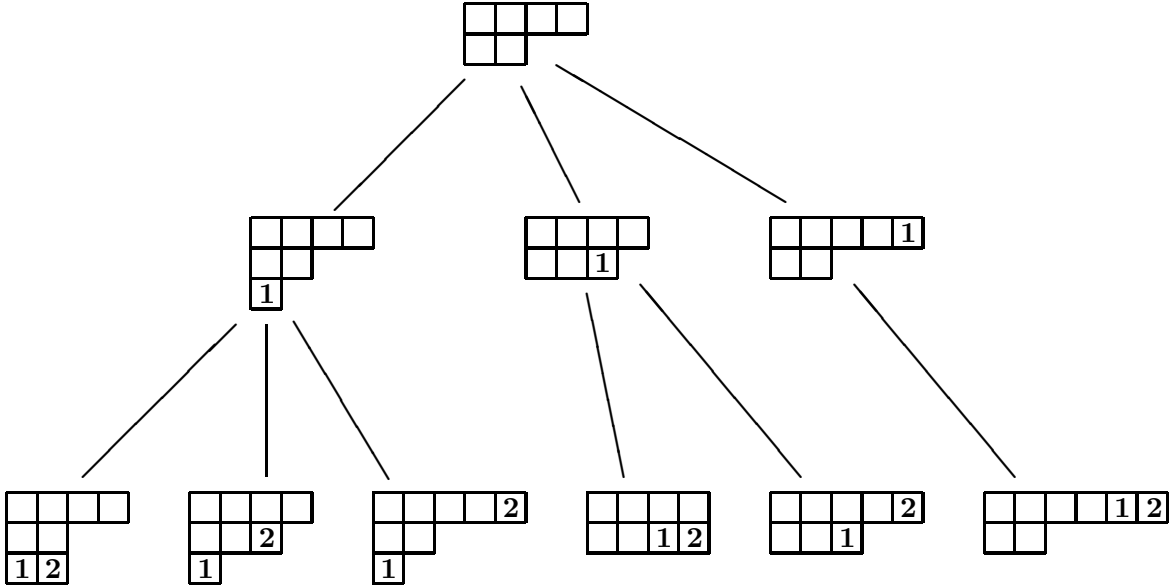

\begin{center}
\sc Appendix B: Application to real enumerative geometry
\end{center}

As in~[So1], explicit rational equivalences can be used 
to find  real solutions to 
problems in enumerative geometry.
Suppose $k={\bf R}$ and $\Fdot$ is a real flag.
\medskip

\noindent{\bf Lemma B.1}
{\em 
Let $\Fdot$ be a flag in ${\bf R}^n$,
 $\alpha \in {[n]\choose m}$, and  $b \geq 0$.
Suppose $X$ is a subvariety of $\,{\bf G}_mV$ which meets 
$Y_{\alpha ,b}\Fdot$ transversally in $d$ points, 
with $e$ real.   
Then there exists $L\in {\bf G}_{n+1-m-b}{\bf R}^n$
with $\Omega_\alpha\Fdot \bigcap \Omega_L$ a generically 
transverse and irreducible intersection such that 
$X$ meets  $\Omega_\alpha\Fdot \bigcap \Omega_L$ transversally in $d$
points, with $e$ real.
}\medskip

\noindent{\bf Proof:} The key idea is that the number of real points in 
the fibre of a finite real morphism is preserved under
real deformations.

We use the notation of \S 3.7.
For $1\leq i\leq b$, let ${\cal G}_{i;{\bf R}}\subset {\cal G}_i$ 
be those cycles $Y_{\alpha,i}(\Fdot,L)$ with $L$ a real 
subspace.
Similarly, let ${\cal G}_{0;{\bf R}}\subset {\cal G}_0$ be the
cycles $\Omega_\alpha\Fdot\bigcap \Omega_L$ with $L$ real.
The description of the maps $\phi_i$ show they are real algebraic.

Suppose  $(\Sigma_i)_0 = Y_{\alpha,i}(\Fdot,M_i)$
meets $X$ transversally in $d$ points, with $e$ real.
For $t$ in ${\bf RP}^1-\{0\}$, $(\Sigma_i)_t\in {\cal G}_{i-1;{\bf R}}$, 
since $\Sigma_i$  is induced by a real 
pencil $L_{i,t}$.
Since $(\Sigma_i)_0$ meets $X$ transversally in $d$ points, 
there is a Zariski open
subset $U$ of ${\bf P}^1$ containing 0 
over which each fibre $(\Sigma_i)_t$ of $\Sigma_{i}$ meets $X$ transversally
in $d$ points.
Let  $I$ be the connected component of $U({\bf R})$ containing 0.
For  $t\in I$, the fibres $(\Sigma_i)_t$ meet $X$ transversally in 
$d$ points, with $e$ real, as
the number of real points 
in the intersection is constant
on connected components of $U({\bf R})$.

Downward induction on $i$ completes the proof,
where, to continue, either reparameterize ${\bf P}^1$ so 
$\infty \in I$, or choose $M_{i-1}$ to be $L_{i,t}$ for some $t\in I-\{0\}$.
\QED\vspace{10pt}

For $\alpha\in{[n]\choose m}$, let $\alpha^\vee\in{[n]\choose m}$ be
$n+1-\alpha_m>\cdots>n+1-\alpha_1$.
Then $\alpha$ and $\alpha^\vee$ are the indices of Poincar\'e dual classes.
\medskip

\noindent{\bf Theorem B.2}
 {\em 
Let $\alpha,\beta \in {[n]\choose m}$ and $a,b,c$ be positive integers 
with 
$a+b+c+|\alpha|+|\beta|\ =\ m(n-m)$,
the dimension of $\,{\bf G}_m{\bf R}^n$.
Let $d$ be the number of pairs 
$(\gamma,\delta)\in (\alpha*a)\times(\beta*b)$
satisfying $\delta^\vee \in \gamma*c$.
Then there exist flags $\Fdot$ and $\Fpdot$ and subspaces $A$, $B$, and 
$C$ of $\,{\bf R}^n$ of respective dimensions $n+1-m-a,n+1-m-b$, and 
$n+1-m-c$ such that 
$\Omega_\alpha\Fdot$, $\Omega_\beta\Fpdot$, 
$\Omega_A$, $\Omega_B$, and $\Omega_C$ meet transversally in $d$ real points.
}\medskip

The conditions on $\alpha,\beta,a,b,c$ ensure that for 
general $\Fdot$, $\Fpdot$, $A$, $B$, and $C$, the 
Schubert varieties meet in a degree $d$ zero cycle.
In fact, we prove this below.
\medskip

\noindent{\bf Proof:}
First suppose  $a=b=0$ so $\Omega_A = \Omega_B = {\bf G}_m{\bf R}^n$.   
This triple intersection arises in the classical proof of Pieri's 
formula ({\em cf.}~[HP]):
Suppose  $\Fdot$, $\Fpdot$,  and  $C$ are in linear general 
position and $c+|\alpha|+|\beta|= m(n-m)$.
If $\beta^\vee\in \alpha*c$, the  varieties 
$\Omega_\alpha\Fdot$, $\Omega_\beta\Fpdot$, and $\Omega_C$
meet transversally in a single
$m$-plane $H_{\alpha\,\beta}$, 
otherwise their intersection is empty.

Let $\beta^\vee\in \alpha*c$.
Then $H_{\alpha\,\beta}$ has a basis $f_1,\ldots,f_m$ with 
$f_j\in K_j := F_{\alpha_j}\cap F'_{\beta_{m+1-j}}$,
where $f_1+\cdots+f_m$ spans the line
$C\bigcap \left( K_1\oplus\cdots\oplus K_m\right)$.
The conditions on $\alpha,\beta$, and $c$ ensure that when 
$\Fdot$,  $\Fpdot$, and $C$ are in linear general position,
each $K_i\neq \{0\}$,  the sum
$K_1+\cdots+K_m$ is direct, the dimension of 
$C\bigcap \left( K_1\oplus\cdots\oplus K_m\right)$
is 1,  that $f_i\not\in F_{\alpha_i+1}$, and that 
$f_i\not\in F'_{\beta_{m+1-i}+1}$.
Thus $H_{\alpha\,\beta}$ is real and if 
$(\alpha,\beta)\neq (\alpha',\beta')$, then 
$H_{\alpha\,\beta}\neq H_{\alpha'\,\beta'}$.
\smallskip

We return to the general case.
The previous discussion shows that if $\Fdot,\Fpdot$, and $C$ are 
real subspaces in linear general position, then 
$$
\left( \sum_{\gamma\in \alpha*a} \Omega_\gamma\Fdot \right)\bigcap
\left( \sum_{\delta\in \beta*b} \Omega_\delta\Fpdot \right)\bigcap
\Omega_C
$$
is a transverse intersection consisting of the $d$ real $m$-planes:
$$
\{H_{\gamma\,\delta}\,|\, \gamma\in \alpha*a,\ \delta\in\beta*b,\
\mbox{ and }\ \delta^\vee \in \gamma*c\}.
$$
By Lemma~B.1, with 
$Y=\left( \sum_{\delta\in \beta*b} \Omega_\delta\Fpdot \right)\bigcap
\Omega_C$, there exists  $A\in {\bf G}_{n+1-m-a}{\bf R}^n$ with 
$$
\Omega_\alpha\Fdot\bigcap\Omega_A\bigcap
\left( \sum_{\delta\in \beta*b} \Omega_\delta\Fpdot \right)\bigcap
\Omega_C
$$
a transverse intersection consisting of $d$ real $m$-planes.
Another application of Lemma~B.1 with 
$Y = \Omega_\alpha\Fdot\bigcap\Omega_A\bigcap\Omega_C$
shows there exists  $B\in  {\bf G}_{n+1-m-a}{\bf R}^n$ with 
$$
\Omega_\alpha\Fdot\bigcap
\Omega_\beta\Fpdot\bigcap \Omega_A\bigcap\Omega_B\bigcap \Omega_C
$$
a transverse intersection consisting of $d$ real $m$-planes.
\QED

\end{document}